\DeclareFontFamily{U}{wncy}{}
\DeclareFontShape{U}{wncy}{m}{n}{%
   <5>wncyr5%
   <6>wncyr6%
   <7>wncyr7%
   <8>wncyr8%
   <9>wncyr9%
   <10>wncyr10%
   <11>wncyr10%
   <12>wncyr6%
   <14>wncyr7%
   <17>wncyr8%
   <20>wncyr10%
   <25>wncyr10}{}
\DeclareMathAlphabet{\cyrille}{U}{wncy}{m}{n}
\def\sh{
\setlength{\unitlength}{.5 pt}
\begin{picture}(40,20)
\put(10,2){\line(1,0){20}}
\put(10,2){\line(0,1){10}}
\put(20,2){\line(0,1){10}}
\put(30,2){\line(0,1){10}}
\end{picture}}

\newcommand{\dlangle}{\langle\!\langle}
\newcommand{\drangle}{\rangle\!\rangle}

\documentclass[a4paper,11pt]{article}

\usepackage{latexsym}
\usepackage{amsmath,amsfonts,amssymb,amsthm}
\newtheorem{thm}{Theorem}[section]
\newtheorem{lem}[thm]{Lemma}

\newtheorem{cor}[thm]{Corollary}
\newtheorem{prop}[thm]{Proposition}

\newtheorem{rem}[thm]{Remark}
\newtheorem{exple}[thm]{Example}

\input epsf

\title{The special subgroup of invertible non-commutative 
rational power series as a metric group}
\author{Roland Bacher
}
\begin{document}
\maketitle

{\sl Abstract: We give a sligthly more natural proof 
of Sch\"utzenberger's Theorem
stating that non-commutative formal power series are rational if and
only if they are recognisable. A byproduct of this proof is
a natural metric on a subgroup of invertible
rational non-commutative power series. We describe a few features 
of this metric group.}

\section{Introduction}

Rational power series in $\mathbb K[[X]]$ over a fixed commutative field
$\mathbb K$ can either be defined
as series representing quotients $\frac{f}{g}$ of two suitable
polynomials $f,g\in\mathbb K[X]$ 
or as ordinary generating series $\sum_{n=0}^\infty s_n X^n$
associated to a sequence satisfying a linear recursion 
relation $s_n=\sum_{j=1}^d \kappa_j s_{n-j}$ for $n\geq N$.
In several non-commuting variables these two descriptions
lead to the notions of rational power series and of 
recognisable power series. Although seemingly distinct, they coincide
for a finite number of variables by a theorem of Sch\"utzenberger.
 
This paper has several goals: Section \ref{sectproofschuetz}
contains an easy proof that rational series are recognisable.
This proof is, up to conventions and notations, the proof 
given in \cite{BR}, except for a slight variation at the end.

This variation consists in an identity which suggests to consider a
natural metric on the multiplicative group of non-commutative
rational power series with constant coefficient $1$. Section 
\ref{sectmetricgroup} describes this metric group.
In particular, we compute the induced metric on the
group generated by $1+X_1,\dots, 1+X_k$ corresponding to 
the image of the Magnus representation of the free group on
$k$ elements. We give also some formulae 
related to the enumeration
of all elements of given norm if the field $\mathbb K$ is finite.

Other parts of this paper discuss enumerative or algorithmic aspects.

The paper is organised as follows:

Section \ref{sectionpows} recalls the basic definitions and 
states Sch\"utzenberger's Theorem.

Section \ref{sectionrecclosure} introduces the notions of recursive
closure and complexity, used as the main tool in the sequel. 
It contains all the necessary ingredients for the proof, given in
Section \ref{sectproofschuetz}, that rational series are recognisable.

Section \ref{sectpresentation} shows how to describe 
``rational'' series using finite amounts of data.
For the sake of completeness, it reproduces also a well-known proof 
of the easy direction of Sch\"utzenberger's Theorem. 

Section \ref{sectproofschuetz} gives an easy proof of 
the ``tedious'' direction of Sch\"utzenberger's Theorem.
It contains also a few formulae useful in the sequel.

Sections \ref{sectionpows}-\ref{sectproofschuetz} contain no 
original results (except perhaps Proposition \ref{propratinverse}
and Corollary \ref{corinvcompl}) and have a large overlap with the first
Chapters of \cite{BR}, except for a few conventions and notations.

Section \ref{secnormal} introduces normal forms. We use them 
for giving some formulae for the number of rational series
of given complexity over finite fields.

Section \ref{sectsaturlevel} addresses a few algorithmic issues.

Section \ref{sectmetricgroup} is devoted to the description and
study of a metric subgroup in the algebra of 
non-commutative rational formal power series. In particular,
we compute this metric on the subgroup defined by the Magnus
representation of a free group. This gives a 
new proof of faithfulness of the Magnus representation of
free groups. At the end of this Section we address enumerative 
questions over finite fields.

The last Section overviews briefly a few related algebraic structures
and recalls a few well-known results concerning linear substitutions,
involutive antiautomorphisms, derivations, Hadamard products, 
shuffle products, compositions, automatic sequences and 
regular languages.

\section{Power series in free non-commuting variables}
\label{sectionpows}

This Section recalls a few basic and well-known facts concerning
formal power series in free non-commuting variables, 
see for instance \cite{St2}, \cite{BR} and \cite{SaSo}.
We try to refer to \cite{BR} and \cite{SaSo} at relevant 
places.
We use sometimes a different terminology, motivated by \cite{B}.

We denote by $\mathcal X^*$ the free monoid over a finite set 
$\mathcal X=\{X_1,X_2,\dots\}$. 
We use boldface capitals 
$\mathbf X,\mathbf T,\mathbf S,\dots$ for 
non-commutative monomials $X_{i_1}X_{i_2}\cdots X_{i_l}
\in\mathcal X^*$. We denote by
$$A=\sum_{\mathbf X\in\mathcal X^*}(A,\mathbf X)\mathbf X$$
a non-commutative formal power series where
$\mathcal X^*\ni\mathbf X\longmapsto (A,\mathbf X)\in \mathbb K$
stands for the coefficient function. The vector space 
$\mathbb K\dlangle \mathcal X\drangle$ consisting of all formal
non-commutative series is an algebra for the convolution
product
$$AB=\sum_{\mathbf X,\mathbf Y\in \mathcal X^*}
(A,\mathbf X)(B,\mathbf Y)\mathbf {XY}$$
of $A,B\in \mathbb K\dlangle \mathcal X\drangle$.

The algebra $\mathbb K\dlangle
\mathcal X\drangle$ contains the subalgebra 
$\mathbb K\langle \mathcal X\rangle$ of non-commutative 
polynomials. The algebra $\mathbb K\langle \mathcal X\rangle$
can also be considered as the the free 
(non-commutative) associative algebra over $\mathcal X$ or as 
the monoid-algebra
$\mathbb K[\mathcal X^*]$ of the free monoid $\mathcal X^*$.

The {\it augmentation map} $\epsilon:\mathbb K\dlangle 
\mathcal X\drangle\longrightarrow\mathbb K$ is the homomorphism of
algebras which sends a series $A$ to its constant coefficient
$\epsilon(A)=(A,\mathbf 1)$. It has a natural section
given by the obvious inclusion $\mathbb K\subset \mathbb K\dlangle
\mathcal X\drangle$ which identifies the field $\mathbb K$
with constant series in $\mathbb K\dlangle \mathcal X\drangle$.
The kernel
$$\ker(\epsilon)=\mathfrak m=\{A\in\mathbb K\dlangle \mathcal X\drangle\ \vert
\ \epsilon(A)=0\}
\subset \mathbb K\dlangle \mathcal X\drangle $$
of $\epsilon$ is the maximal ideal consisting of all 
formal power series without
constant coefficient of the local algebra $\mathbb K\dlangle 
\mathcal X\drangle$. 

$\mathbb K\dlangle\mathcal X\drangle$ is a complete topological space
for the neighbourhood filter $A+\mathfrak m^i,\ i=0,1,2,\dots$
of 
$A\in\mathbb K\dlangle\mathcal X\drangle$.

We have for $a\in\mathfrak m$ the equality
$$(1-a)(1+\sum_{n=1}^\infty a^n)=1\ .$$
It shows that a formal power series 
$A\in\mathbb K\dlangle \mathcal X\drangle $ is invertible 
with respect to the non-commutative product if and only if
$A\not\in \mathfrak m$.

We denote by $\mathbb K\dlangle \mathcal X\drangle^*=
\epsilon^{-1}(\mathbb K^*)=\mathbb K^*+\mathfrak m$ the 
non-commutative {\it group of units} of the
algebra $\mathbb K\dlangle \mathcal X\drangle $
formed by all invertible elements. 
We call the subgroup $S\mathbb K\dlangle \mathcal X\drangle^*=\{A\in 
\mathbb K \dlangle \mathcal X\drangle^*\ \vert\epsilon(A)=1\}$
the {\it special group of units}. The homomorphism
$$\mathbb K\dlangle\mathcal X\drangle^*\ni A\longmapsto
\frac{1}{\epsilon(A)}A\in S\mathbb K\dlangle\mathcal X\drangle^*$$
identifies $S\mathbb K\dlangle\mathcal X\drangle^*$ with the 
projective quotient-group
$\mathbb K\dlangle\mathcal X\drangle^*/\mathbb K^*$
and shows the direct product decomposition
$\mathbb K\dlangle\mathcal X\drangle^*=\mathbb K^*\times S\mathbb K\dlangle
\mathcal X\drangle^*$. 

\begin{prop} \label{proptorsioninKX}
The roots of $1$ contained in the central subgroup $\mathbb K^*$
of $\mathbb K\dlangle\mathcal X\drangle^*$ are the only torsion elements 
of $\mathbb K\dlangle\mathcal X\drangle^*$.
\end{prop}

\begin{cor}
The group $S\mathbb K\dlangle\mathcal X\drangle^*$ is without torsion.
\end{cor}

{\bf Proof of Proposition \ref{proptorsioninKX}} 
Let $\alpha(1+a)\in\mathbb K\dlangle\mathcal X\drangle^*$
be a torsion element of order exactly 
$d$ with $\alpha\in\mathbb K^*$ and $a\in
\mathfrak m$. This implies $\alpha^d=1$ and $(1+a)^d=1$.
It is thus enough to show that we have $a=0$.
We have 
$$1=(1+a)^d=1+\sum_{j=1}^d{d\choose j}a^j$$
which implies $a=0$ by considering in positive characteristic
the smallest strictly positive
integer $j'$ such that ${d\choose j'}$ is not divisible by the 
characteristic of $\mathbb K$.\hfill$\Box$

\begin{rem} Bourbaki, see for example Page 45 of \cite{Bour23},
calls $\mathbb K\dlangle\mathcal X\drangle$ the Magnus algebra 
and the group of units $\mathbb K\dlangle\mathcal X\drangle^*$ the
Magnus group. We do not follows this terminology.
\end{rem}

\begin{rem} The theory of power series in non-commuting variables
can be developped over an associative semi-ring, cf. \cite{BR}
and \cite{SaSo}. 
From the point of
view of the associated unit group there is no loss of generality
by requiring all coefficients (except perhaps the constant
coefficient) to belong to the maximal subring of such a semi-ring.
One can indeed show that a series of the form $1+a$ with $a\in
\mathfrak m$ is invertible if and only if all coefficients of 
$a$ have additive inverses.
\end{rem}

\subsection{Rational series}


We use the convention that algebras are with unit. 
In particular, every subalgebra
of $\mathbb K\dlangle \mathcal X\drangle$ contains 
the field $\mathbb K$.

A subalgebra of $\mathbb K\dlangle\mathcal X\drangle$ 
is {\it rationally closed} or {\it full} if
it intersects the unit group $\mathbb K\dlangle\mathcal X
\drangle^*$ in a subgroup. 
The {\it rational closure} of a subset $\mathcal S\subset \mathbb
K\dlangle \mathcal X\drangle$ is the smallest rationally 
closed subalgebra of 
$\mathbb K\dlangle \mathcal X\drangle$ which contains 
$\mathcal S$. 
 
The rational closure $\mathbb K\dlangle \mathcal X\drangle_{rat}$ 
of $\mathcal X$ is called the {\it algebra of rational
  series} or the {\it rational subalgebra} of $\mathbb K\dlangle
\mathcal X\drangle$ and is formed by {\it rational elements}. 
It is the smallest subalgebra of $\mathbb K\dlangle \mathcal X\drangle$
which contains the polynomial subalgebra 
$\mathbb K\langle \mathcal X\rangle$ and intersects the 
unit group $\mathbb K\dlangle \mathcal X\drangle^*$ in a subgroup
$\mathbb K\dlangle \mathcal X\drangle^*_{rat}$, called
the {\it group of rational units}. We denote by
$S\mathbb K\dlangle\mathcal X\drangle^*_{rat}=
S\mathbb K\dlangle \mathcal X\drangle^*\cap 
\mathbb K\dlangle\mathcal X\drangle_{rat}$ the
{\it special group of rational units}. We write
$S\mathbb K\dlangle \mathcal X\drangle^*_{pol}$ for
the subgroup of $S\mathbb K\dlangle \mathcal X\drangle^*_{rat}$
generated by all elements in $S\mathbb K\dlangle \mathcal X\drangle_{rat}^*
\cap \mathbb K\langle \mathcal X\rangle$. It follows for example from
Chapter IV, Section 3 of \cite{BR} that 
$S\mathbb K\dlangle \mathcal X\drangle^*_{pol}$ is 
a proper subgroup of $S\mathbb K\dlangle \mathcal X\drangle_{rat}^*$
if $\mathcal X$ contains more than one variable.

\begin{rem} If $\mathcal X$ is reduced to a unique element $X$,
the algebra $\mathbb K\dlangle X\drangle$ is the commutative
algebra $\mathbb K[[X]]$ of ordinary formal power series in one 
variable and we have $S\mathbb K\dlangle X\drangle_{pol}^*=
S\mathbb K\dlangle X\drangle_{rat}^*$. 
\end{rem}


\begin{rem} The groups $S\mathbb R\dlangle\mathcal X\drangle^*$ and 
$S\mathbb C\dlangle\mathcal X\drangle^*$ are infinite-dimensional 
real or complex Lie-groups
with Lie algebra $\mathfrak m$ and Lie-bracket $[a,b]=ab-ba$.
\end{rem}

\subsection{Recognisability and Sch\"utzenberger's Theorem}

An element $A\in\mathbb K\dlangle \mathcal X\drangle$ is
{\it recognisable} if there exists a finite-dimensional 
$\mathbb K-$vector space $\mathcal V$, a morphism of monoids
$\mu:\mathcal X^*\longrightarrow \hbox{End}(\mathcal V)$
and elements $\alpha\in \mathcal V,\ \omega\in\hbox{Hom}(\mathcal V,
\mathbb K)$ such that
$$A=\sum_{\mathbf X\in\mathcal X^*}
\omega(\mu(\mathbf X)\alpha)\ \mathbf X\ .$$

The following result is due to Sch\"utzenberger, see for example
Theorem 6.5.7 in \cite{St2}, Theorem 7.1, Page 15 in \cite{BR} or
Theorem 2.3, Page 22 in \cite{SaSo}.

\begin{thm}\label{thmSchutz} 
Given a finite set $\mathcal X$ of
free non-commuting variables, an element
$A\in\mathbb K\dlangle\mathcal X\drangle$ is recognisable if and only
if it is rational.
\end{thm}

\begin{rem} Theorem \ref{thmSchutz} does not hold if $\mathcal X$
is an infinite set: A series of the form $A=\sum_{X_j\in\mathcal X}
\lambda_j X_j$ is recognisable. It is however 
not rational if infinitely many coefficients $\lambda_j$ are 
non-zero.

Theorem \ref{thmSchutz} remains however true when considering
only elements of the subalgebra 
$$\mathbb K_f\dlangle \mathcal X\drangle=\bigcup_{\mathcal 
X_f\hbox{ finite subset of }\mathcal X}\ 
\mathbb K\dlangle\mathcal X_f\drangle$$
formed by elements of $\mathbb K\dlangle\mathcal X\drangle$
involving only finitely many variables of $\mathcal X$. 

Such subtleties can be avoided by requiring
finiteness of the set $\mathcal X$ of variables. 
\end{rem}




Since rational elements, recognisable elements and elements of
finite complexity in $\mathbb K\dlangle \mathcal X\drangle$
(for $\mathcal X$ finite)
coincide by Sch\"utzenberger's Theorem, we drop these distinctions 
after completion of the proof of Theorem \ref{thmSchutz}
and speak simply of rational elements.

\section{Recursive closure in $\mathbb K\dlangle \mathcal X\drangle$
and complexity}\label{sectionrecclosure}

This section introduces the notion of recursive closure and identifies
the set of series having a finite-dimensional recursive closure with
the set of recognisable series.

One should mention
that the definition of the recursive closure is not
completely canonical: there are three natural 
choices due to the fact that one can consider shift maps acting 
on the ``right'', on the ``left'' or on 
``both sides''. The differences between right and left are
minor and lead to isomorphic theories (our conventions coincide
with the choice of \cite{SaSo}, the book \cite{BR} uses
the opposite conventions). The theory for
the symmetric choice of a bilateral action
is more cumbersome: ``breaking the symmetry'' makes life
easier. 

\subsection{Recursively closed subspaces and complexity}

The (generalised) {\it Hankel matrix} $H=H(A)$ of a series
$$A=\sum_{\mathbf X\in\mathcal X^*}(A,\mathbf X)\mathbf X\in
\mathbb K\dlangle \mathcal X\drangle $$ 
is the infinite matrix with rows and columns 
indexed by all elements of the free monoid $\mathcal X^*$,
whose entries are given by 
$H_{\mathbf X,\mathbf X'}=(A,\mathbf X\mathbf X')$. 

We associate to the row of index $\mathbf T$ in $H(A)$ the series 
$$\rho(\mathbf T)A=\sum_{\mathbf X\in\mathcal X^*} 
(A,\mathbf{XT}) \mathbf X=\sum_{\mathbf X\in\mathcal X^*} 
H_{\mathbf X,\mathbf T} \mathbf X\in\mathbb K\dlangle \mathcal X\drangle\ .$$
Using the terminology of \cite{B}, we call the vector-space 
$\overline A\subset \mathbb K\dlangle
\mathcal X\drangle$ spanned by $\rho(\mathbf T)A,\ 
\mathbf T\in\mathcal X^*$, the {\it recursive closure} of $A$
and its dimension $\dim(\overline A)\in\mathbb N\cup
\{\infty\}$ the {\it rank} or {\it complexity} of $A$. The complexity 
$\dim(\overline A)$ of a series $A$ can be thought of as a sort 
of ``degree'' of a non-commutative rational series and
is equal to the rank of the Hankel matrix $H(A)$, defined as the
dimension of the vector space spanned by all rows (or, equivalently,
by all columns) of $H(A)$. Let us add that \cite{BR} and \cite{SaSo}
use the terminology ``rank'' instead of complexity. We prefer complexity
in order to avoid confusions related to the matrix-context described
in \cite{B}.

A subspace
$\mathcal V\subset \mathbb K\dlangle\mathcal X\drangle$ is 
{\it recursively closed} if it contains the recursive closure 
of all its elements.

The {\it shift map} of a monomial $\mathbf T\in\mathcal X^*$ 
is the $\mathbb K-$linear map 
$\rho(\mathbf T)\in\mathop{End}(\mathbb K\dlangle\mathcal X\drangle)$
defined as above by $A\longmapsto \rho(\mathbf T)A=
\sum_{\mathbf X\in\mathcal X^*}(A,\mathbf{XT})\mathbf X$.
The identity 
$$\rho(\mathbf T)(\rho(\mathbf{T'})A)=\sum_{\mathbf X\in\mathcal X^*}
(A,\mathbf{XTT'})\mathbf X=\rho(\mathbf{TT'})A$$ 
shows that the shift maps $\rho:\mathcal X^*
\longrightarrow \mathop{End}(\mathbb K\dlangle X\drangle)$
define a linear representation of the free monoid $\mathcal X^*$.
We call this linear representation the {\it shift monoid}.
Since a recursively closed subspace $\mathcal V\subset\mathbb K\dlangle 
\mathcal X\drangle$ is stable under the action of the shift monoid,
restriction of the shift-monoid to $\mathcal V$ yields a 
subrepresentation $\rho_{\mathcal V}:\mathcal X^*\longrightarrow
\mathop{End}(\mathcal V)$, called the shift monoid of $\mathcal V$.
If $\mathcal V=\overline A$ is the recursive closure of an
element $A\in\mathbb K\dlangle \mathcal X\drangle$, we speak simply
of the shift-monoid $\rho_{\overline A}$ of $A$.

\begin{exple}\label{expleunsurunmoinsxy} 
For $A=1/(1-XY)$ we have $\rho(X)A=0$, 
$\rho(Y)A=AX$, $\rho(X)(AX)=A$ and $\rho(Y)(AX)=0$.
The series $A=1+AXY$ is thus of complexity $2$ and has recursive closure
$\overline A=\mathbb K A+\mathbb K AX$. The shift monoid 
$\rho_{\overline A}(\mathcal X^*)$ is generated by the two matrices 
$$\rho_A(X)=\left(\begin{array}{cc}0&1\\0&0\end{array}\right)\hbox{ and }
\rho_A(Y)=\left(\begin{array}{cc}0&0\\1&0\end{array}\right)$$
acting by left multiplication on column vectors $\left(\begin{array}{c}
\alpha\\ \beta\end{array}\right)$ corresponding to $\alpha A+\beta AX\in
\overline A$. We leave it to the reader to check that 
the shift-monoid $\rho_{\overline A}(\mathcal X^*)$ is the finite monoid
consisting of the identity $\rho(\emptyset)$, four non-zero elements 
$\rho_A(X),\rho_A(Y),\rho_A(YX),\rho_A(XY)$ and of the zero element
$\rho_A(XX)=\rho_A(YY)$.

\end{exple}

\begin{exple} \label{remdimratfrac} The rational elements
of the commutative algebra $\mathbb K\dlangle X\drangle=\mathbb K[[X]]$ 
of formal power series in one variable are given by rational fractions
$f/g$ with $g\in\mathbb K\dlangle X\drangle^*$ invertible.

The complexity $\dim(\overline A)$ of a non-zero 
rational fraction $A\in\mathbb K\dlangle X\drangle_{rat}$ in one 
variable equals $\dim(\overline
A)=\mathop{max}(1+\mathop{deg}(f),\mathop{deg}(g))$ where 
$f/g=A$ is a reduced expression for $A$, cf. Exercise 3, Page 60 of 
\cite{SaSo}.

The action of the shift map $\rho(X)$ on $A\in\mathbb K\dlangle X\drangle$
is given by
$$A=\sum_{n=0}^\infty \alpha_nX^n\longmapsto
\rho(X)A=\sum_{n=1}^\infty \alpha_nX^{n-1}$$
and corresponds thus to the well-known unilateral shift
$$(\alpha_0,\alpha_1,\alpha_2,\dots)\longmapsto
(\alpha_1,\alpha_2,\alpha_3,\dots)$$
on the sequence of coefficients of a formal power-series.
\end{exple}

Recognisable series can be characterised by the following result,
cf Proposition 5.1, Page 9, of \cite{BR}.

\begin{prop} \label{proprecog=fincompl}
An element $A\in \mathbb K\dlangle \mathcal X\drangle$ 
is recognisable if and only if it is of finite complexity.
\end{prop}

Sch\"utzenberger's Theorem amounts thus to the assertion that
a series in $\mathbb K\dlangle\mathcal X\drangle$ is rational 
if and only if it has finite complexity.

{\bf Proof of Proposition \ref{proprecog=fincompl}}
The identity 
$$A=\sum_{\mathbf X\in\mathcal X^*}\epsilon(\rho(\mathbf X)A)\mathbf X$$
for $A\in\mathbb K\dlangle\mathcal X\drangle$ implies that an element
of finite complexity is recognisable by considering 
$\mathcal V=\overline A$,
$\mu_A=\rho\vert_{\overline A}:\mathcal X^*\longrightarrow
\hbox{End}(\mathcal V)$, $\alpha=A\in\mathcal V$ and $\omega=
\epsilon\vert_{\overline A}\in \hbox{Hom}(\mathcal V,\mathbb K)$.

On the other hand, consider a recognisable element $A$ given by
$$A=\sum_{\mathbf X\in\mathcal X^*}\omega(\mu_A(\mathbf X)\alpha)\mathbf X
$$
where $\mu_A:\mathcal X^*\longrightarrow\mathrm{End}(\mathcal V)$
is a linear representation of $\mathcal X^*$ on some finite-dimensional
vector space $\mathcal V$ and where $\alpha\in\mathcal V,\ \omega\in
\hbox{Hom}(\mathcal V,\mathbb K)$. The obvious identities 
$$\rho(\mathbf T)A=\sum_{\mathbf X\in \mathcal X^*}
\omega(\mu_A(\mathbf{XT})\alpha)\mathbf X=\sum_{\mathbf X\in \mathcal X^*}
\omega\left(\mu_A(\mathbf{X})(\mu_A(\mathbf T)\alpha)\right)\mathbf X$$
show the inclusion
$$\overline A\subset \{\sum_{\mathbf X\in\mathcal X^*}\omega(\mu_A(
\mathbf X)\beta)\mathbf X\ \vert\ \beta\in\mathcal V\}$$
which implies $\dim(\overline A)\leq \dim(\mathcal V)<\infty$.
\hfill$\Box$ 

\begin{rem} The linear representation $\rho_{\mathcal A}:
\mathcal X^*\longrightarrow
\mathop{End}(\mathcal A)$ associated to a recursively closed 
subspace $\mathcal A\subset \mathbb K\dlangle
\mathcal X\drangle$ extends to the monoid algebra
$\mathbb K[\mathcal X^*]=\mathbb K\langle \mathcal X\rangle$ and is thus,
up to conjugation by an element of $\hbox{Aut}(\mathcal A)$,
uniquely defined by the two-sided 
ideal 
$$\mathcal I_{\mathcal A}=\ker(\rho:\mathbb K\langle \mathcal X\rangle
\longrightarrow \hbox{End}(\mathcal A))\subset
\mathbb K\langle \mathcal X\rangle$$
called the {\it syntaxic ideal} in \cite{BR}. The quotient algebra
$\mathbb K\langle \mathcal X\rangle/\mathcal I_{\mathcal A}$
can be identified with the monoid-algebra $\mathbb K[\rho_{\mathcal A}
(\mathcal X^*)]$. It is called the syntaxic algebra of $A$ in \cite{BR} 
if $\mathcal A=\overline A$. Let me also mention that
the complexity corresponds to the rank of an element in \cite{BR} where
elements of the maximal ideal $\mathfrak m\subset \mathbb K
\dlangle\mathcal X\drangle$ are called proper elements.
\end{rem}

\section{Presentations}\label{sectpresentation}

This section introduces recursive presentations (corresponding to 
linear representations in \cite{BR}) for series of finite complexity. 
It contains no new results and is mainly included for the convenience
of the reader.

A {\it recursive presentation} is a finite system of equations
of the form 
$$\left\lbrace\begin{array}{l}
\displaystyle A_1=\gamma_1+\sum_{i=1}^a A_i\alpha_{i,1},\\
\displaystyle \qquad\vdots\\
\displaystyle A_a=\gamma_a+\sum_{i=1}^a A_i\alpha_{i,a},
\end{array}\right.$$
with unknowns $A_1,\dots,A_a$,
constants $\gamma_1,\dots,\gamma_a\in\mathbb K$ and
homogeneous linear forms
$\alpha_{i,j}\in\mathbb K\langle \mathcal X\rangle\cap
\mathfrak m$ in the variables $\mathcal X$
for $(i,j)\in\{1,\dots,a\}^2$.

\begin{prop} \label{proprecpr=fincompl}
(i) A recursive presentation involving $a$ equations in
$a$ unknowns $A_1,\dots,A_a$ has a unique solution 
$(A_1,\dots,A_a)\in\left(\mathbb K\dlangle \mathcal X\drangle\right)^a$.

\ \ (ii) The series $A_1,\dots,A_a$ defined by the solution of a 
recursive presentation span a recursively closed vector space.
\end{prop}

{\bf Proof} The proof of assertion (i) is by ``bootstrapping'': The inclusions
$\alpha_{i,j}\subset \mathfrak m$ imply that the 
dynamical system of $(\mathbb K\dlangle\mathcal X\drangle)^a$ given by
the map 
$$(\tilde A_1,\dots,\tilde A_j,\dots,\tilde A_a)\longmapsto
(\dots,\gamma_j+\sum_{i=1}^a\tilde A_i\alpha_{i,j},\dots)$$ 
has a unique fixpoint which coincides thus with the 
solution $(A_1,\dots,A_a)$, determined by the recursive presentation.
This fixpoint is attracting for the topology defined by the neighbourhood
filter $\mathfrak m^i,\ i=0,1,2,\dots$ of $0$.

Since the $\alpha_{i,j}$'s are homogeneous linear forms of 
$\mathbb K\langle\mathcal X\rangle$, we have
$$\rho(X)A_j=\sum_{i=1}^a\rho(X)\left(A_i\alpha_{i,j}\right)=
\sum_{i=1}^aA_i\rho(X)\alpha_{i,j}\in\sum_{i=1}^a\mathbb K A_i$$ 
for all $X\in\mathcal X$. This shows $\dim(\overline{A_j})
\leq a$ and ends the proof.\hfill$\Box$

A recursive presentation is {\it reduced} if the series $A_1,\dots,A_a$ 
defined by its solution are linearly independent.
A recursive presentation
with solution $(A_1,\dots,A_a)$ is a {\it recursive presentation of 
$A=A_1$}. A recursive presentation of $A$ is {\it minimal}
if the series $A_1,\dots,A_a$ defined by its solution form a basis 
of $\overline A$.
A minimal recursive presentation of $0\in
\mathbb K\dlangle \mathcal X\drangle$ is by convention
the empty recursive presentation with zero equations and unknowns.

\begin{prop} \label{propfincomplhasrecpres}
Every element $A\in\mathbb K\dlangle\mathcal X\drangle$
with finite complexity has a minimal recursive presentation.
\end{prop}

{\bf Proof} If $A\not=0$, we can complete $A$ to
a basis $A_1=A,\dots,A_a$ of $\overline A$. The result follows now from
the observation that
the equations
$$\begin{array}{l}
\displaystyle A_1=\epsilon(A_1)+\sum_{X\in\mathcal X}\left(\rho(X)A_1
\right)X,\\
\displaystyle \qquad \vdots\\
\displaystyle A_a=\epsilon(A_a)+\sum_{X\in\mathcal X}\left(\rho(X)A_a
\right)X\end{array}$$
define a minimal recursive presentation of $A_1=A$.\hfill$\Box$

\begin{rem} The set of minimal presentations of a non-zero element
$A$ with finite complexity is in bijection with the set of sequences
$A_1=A,A_2,\dots\subset \overline A$ extending $A_1=A$ to a 
basis of $\overline A$.
\end{rem}

\begin{exple} \label{explencXY}
Setting $A_1=A=1/(1-XY)$ and $A_2=\rho(Y)A=AX$,
the rational series $A=A_1=1/(1-XY)$ of Example \ref{expleunsurunmoinsxy}
is defined by the minimal recursive presentation
$$A_1=1+A_2Y,\qquad A_2=A_1X\ .$$
\end{exple}

The following result is the ``easy'' direction of Sch\"utzenberger's
Theorem.

\begin{prop} \label{proprecognimpliesrational}
A recognisable series of $\mathbb K\dlangle\mathcal X\drangle$
is rational.
\end{prop}

\begin{exple} The recursive presentation
$$\begin{array}{l}
\displaystyle A=1+BX+(A+B)Y,\\
\displaystyle B=1+(A+B)X+AY\end{array}$$
implies
$$\begin{array}{l}
A=1+X+2Y+2X^2+YX+3XY+3Y^2+3X^3+\dots\\
B=1+2X+Y+3X^2+3YX+XY+2Y^2+5X^3+\dots\end{array}$$
and defines by Propositions \ref{proprecpr=fincompl} 
and \ref{proprecog=fincompl}
a recognisable series $A\in \mathbb K\dlangle X,Y\drangle$
which is rational by Proposition \ref{proprecognimpliesrational}.
Eliminating $B$ in the recursive presentation given above
yields indeed the rational expression
$$A=\left(1+\frac{1}{1-X}(X+Y)\right)\left(1-Y-(X+Y)\frac{1}{1-X}(X+Y)
\right)^{-1}\ .$$
\end{exple}

{\bf Proof of Proposition \ref{proprecognimpliesrational}}
The following proof by non-commutative Gaussian elimination 
is borrowed from \cite{St2}.

A recognisable series $A\in\mathbb K\dlangle \mathcal X\drangle$ 
is of finite complexity by Proposition \ref{proprecog=fincompl}.
It is thus defined by a recursive presentation by
Proposition \ref{propfincomplhasrecpres}. We can thus suppose that 
$A=A_1$ is given by a system of equations of the form
$$\left\lbrace\begin{array}{l}
\displaystyle A_1=\gamma_1+\sum_{i=1}^a A_i\alpha_{i,1},\\
\displaystyle \qquad\vdots\\
\displaystyle A_a=\gamma_j+\sum_{i=1}^a A_i\alpha_{i,a},
\end{array}\right.$$
with $\gamma_1,\dots,\gamma_a\in\mathbb K\dlangle
\mathcal X\drangle_{rat}$ and
$\alpha_{i,j}\in\mathbb K\dlangle \mathcal X\drangle_{rat}\cap
\mathfrak m$ for $i,j\in\{1,\dots,a\}$ since a presentation is a particular
case of such a system of equations.
Solving the last equation for $A_j$ we get
$$A_a=\left(\gamma_a+\sum_{i=1}^{a-1}A_i\alpha_{i,a}\right)\frac{1}{1-
\alpha_{a,a}}$$ where the assumption $\alpha_{a,a}\in\mathbb K
\dlangle \mathcal X\drangle_{rat}\cap \mathfrak m$ implies 
$1/(1-\alpha_{j,j})\in\mathbb K\dlangle \mathcal X\drangle_{rat}$.
If $a=1$ we have $A_1=A_a=\gamma_1/(1-\alpha_{1,1})
\in\mathbb K\dlangle\mathcal X\drangle_{rat}$ and we are done.
Otherwise, we get by elimination of $A_a$ the system
of equations 
$$\left\lbrace\begin{array}{l}
\displaystyle A_1=\tilde \gamma_1+\sum_{i=1}^{a-1} A_i\tilde\alpha_{i,1},\\
\displaystyle \qquad\vdots\\
\displaystyle A_{a-1}=\tilde\gamma_j+\sum_{i=1}^{a-1} A_i\tilde\alpha_{i,a-1},
\end{array}\right.$$
with
$$\tilde \gamma_j=\gamma_j+\gamma_a\frac{1}{1-\alpha_{a,a}}\alpha_{a,j}
\in\mathbb K\dlangle \mathcal X\drangle_{rat}$$
and
$$\tilde \alpha_{i,j}=\alpha_{i,j}+\alpha_{i,a}\frac{1}{1-\alpha_{a,a}}
\alpha_{a,j}\in\mathbb K\dlangle\mathcal X\drangle_{rat}\cap\mathfrak m$$
for $i,j\in\{1,\dots,a-1\}$. This proves the result by induction on the
number of equations and unknowns.\hfill$\Box$

\subsection{Reducing recursive presentations}\label{secreducrec}

A recursive presentation 
$$A_j=\gamma_j+\sum_{i=1}^aA_i\alpha_{i,j},\ j=1,\dots,a$$
with solution $A_1,\dots,A_a\in\mathbb K\dlangle\mathcal X\drangle$
is reduced if $\mathcal K=\{0\}$ where $\mathcal K
\subset \mathbb K^a$ is the kernel of the map
$\pi:\mathbb K^a\longrightarrow \mathbb K\dlangle\mathcal X\drangle$
defined by
$$\mathbb K^a\ni\lambda=(\lambda_1,\dots,\lambda_a)\longmapsto
\pi(\lambda)=\sum_{j=1}^a\lambda_jA_j\in
\mathbb K\dlangle\mathcal X\drangle\ .$$
This kernel $\mathcal K$ can be computed as follows: Let $\mathcal K_0$
denote the kernel of the linear form
$$\lambda=(\lambda_1,\dots,\lambda_a)\longmapsto \epsilon\circ
\pi(\lambda)=\epsilon\left(\sum_{j=1}^a\lambda_j A_j\right)=
\sum_{j=1}^a\lambda_j \gamma_j$$
corresponding to the image of $\sum_{j=1}^a\mathbb KA_j$ under
the augmentation map $\epsilon:\mathbb K\dlangle\mathcal X\drangle
\longrightarrow \mathbb K\dlangle\mathcal X\drangle/\mathfrak m$. 

One defines now recursively
$\mathcal K_{i+1}\subset \mathcal K_i$ as the intersection
$$\mathcal K_{i+1}=\mathcal K_i\cap \bigcap_{j=1}^a\tilde\rho(X_j)^{-1}
(\mathcal K_i)$$
where $\tilde\rho(X_j):\mathbb K^a\longmapsto \mathbb K^a$
is the linear application defined by 
$$e_h\longmapsto \sum_{i=1}^a (X_j,\alpha_{i,h})e_i$$
with respect to the standard basis $e_1,\dots,e_a$ of 
$\mathbb K^a$.
Since $\mathcal K_0$ is finite-dimensional, the sequence 
$$\mathcal K_0\supset \mathcal K_1\supset \mathcal K_2\supset
\dots$$ stabilises and the definition of $\mathcal K_{i+1}$
shows that $\mathcal K_{h}=\mathcal K_{h+1}$ implies $\mathcal K_i=
\mathcal K_h$ for all $i\geq h$. We set $\mathcal K_\infty=
\mathcal K_h$ for such an integer $h$.

\begin{prop} The application 
$$\lambda=(\lambda_1,\dots,\lambda_a)
\longmapsto \pi(\lambda)=\sum_{j=1}^a\lambda_j A_j$$
defines an isomorphism from $\mathbb K^a/\mathcal K_\infty$
onto the recursively closed vector space 
$\sum_{j=1}^a\mathbb K A_j\subset \mathbb K\dlangle\mathcal X\drangle$.
\end{prop}

{\bf Proof} The maps $\tilde\rho(\mathcal X)\in\hbox{End}(\mathbb K^a)$
define sections of the shift maps $\rho(\mathcal X)\in
\hbox{End}(\sum_{j=1}^\infty \mathbb KA_j)$. Since $A\in\mathbb K\dlangle
\mathcal X\drangle$ is zero if and only if $\epsilon(\rho(\mathbf X)A)=0$
for all $\mathbf X\in\mathcal X^*$, the result follows from the observation
that $\mathcal K_\infty$ is the largest $\tilde\rho(\mathcal X)-$stable
subspace of $\pi^{-1}(\mathfrak m)=\ker(\epsilon\circ \pi)$.
\hfill$\Box$ 

\begin{exple} For the recursive presentation 
$$\left\lbrace\begin{array}{l}
\displaystyle A_1=1+(A_1+A_2)X+(A_1-A_3)Y\\
\displaystyle A_2=-1+A_2X+(A_3-A_2)Y\\
\displaystyle A_3=(A_2+A_3)X+(A_1-A_2)Y
\end{array}\right.$$
we get $\mathcal K_0=\mathbb K(A_1+A_2)+\mathbb K A_3$.
The computations 
$$\epsilon\left(\tilde\rho(X)(A_1+A_2)\right)=
\epsilon(A_1+2A_2)=-1$$
and 
$$\epsilon(\tilde\rho(X)A_3)=\epsilon(A_2+A_3)=-1$$
show $\mathcal K_1\subset \mathbb K(A_1+A_2-A_3)$ and
we have
$$\tilde\rho(X)(A_1+A_2-A_3)=A_1+A_2-A_3,\ \tilde\rho(Y)
(A_1+A_2-A_3)=0$$
implying $\mathcal K_\infty =\mathbb K(A_1+A_2-A_3)$.
We have thus the relation $A_3=A_1+A_2$ which we can use
to obtain the reduced recursive presentation 
$$\left\lbrace\begin{array}{l}
\displaystyle A_1=A+(A_1+A_2)X-A_2Y\\
\displaystyle A_2=-1+A_2X+A_1Y\end{array}\right.$$
of the recursively closed vector-space 
$\mathbb KA_1\oplus \mathbb KA_2=\sum_{j=1}^3\mathbb KA_j
\subset \mathbb K\dlangle X,Y\drangle$.
\end{exple}

\subsection{Minimal presentations}\label{secminimal}

A reduced recursive presentation (defining 
$a$ linearly independent series $A_1,\dots,A_a$) of a series $A=A_1$ is 
not necessarily minimal since the inclusion $\overline A\subset 
\oplus_{j=1}^a \mathbb K A_j$ can be strict. 

A minimal presentation of such a series $A\in
\mathbb K\dlangle X_1,\dots,X_k\drangle$ can be constructed as follows:
Set $\mathcal A_0=\mathbb KA$ and define $\mathcal A_{i+1}$
recursively by
$$\mathcal A_{i+1}=\mathcal A_i+\sum_{j=1}^k 
\mathbb K\rho(X_j)\mathcal A_i\ .$$ 
The inclusion $\mathcal A_{i+1}\subset \oplus_{j=1}^a
\mathbb KA_j$ shows that 
there exists an integer $h$ such that $\mathcal A_h=
\mathcal A_{h+1}$. We set $\mathcal A_\infty=\mathcal A_h$ since 
$\mathcal A_i=\mathcal A_h$ for all $i\geq h$. 

We have the following obvious result, given without proof:

\begin{prop} We have $\overline A=\mathcal A_\infty$.
\end{prop}

It is now easy to construct a minimal presentation of $A$ by
extending $\tilde A_1=A$ (for $A\not=0$)
to a basis $\tilde A_1=A,\tilde A_2,\dots$ of $\mathcal A_\infty$.

\section{Proof of Sch\"utzenberger's Theorem}\label{sectproofschuetz}

This Section contains a proof along the lines of \cite{BR}
except fo a small variation involving Proposition \ref{propratinverse}
of the fact that rational
series are recognisable. The idea is to use Proposition 
\ref{proprecog=fincompl} and to show 
that the set of series of finite complexity is a rationally
closed subalgebra of $\mathbb K\dlangle\mathcal X\drangle$.
The direction ``rational implies recognisable'' of
Sch\"utzenberger's Theorem follows then from the obvious remark that
non-commutative polynomials have finite complexity.

Given subspaces $E,F\subset \mathbb K\dlangle \mathcal X\drangle$,
we denote by
$$\begin{array}{l}
\displaystyle E+F=\{A+B\ \vert A\in E,B\in F\},\\
\displaystyle EF=\{\sum_{j}A_jB_j\ \vert A_j\in E,B_j\in F\}
\end{array}$$
the subspaces of $\mathbb K\dlangle \mathcal X\drangle$
spanned by sums, respectively products, of an element
in $E$ and an element in $F$. 

The following lemma, corresponding to Lemme 7.2, Page 15 of \cite{BR},
is a key ingredient.

\begin{lem} \label{lemrhoXAB} (i) We have
$$\rho(X)(AB)=\epsilon(B)\rho(X)A+A(\rho(X)B)$$
for all $X\in\mathcal X$.

\ \ (ii) For $A\in\mathbb K\dlangle\mathcal X\drangle^*$
invertible and for $X\in\mathcal X$, we have 
$$\rho(X)(A^{-1}B)=A^{-1}\left(-\epsilon(B)/\epsilon(A)\rho(X)A+
\rho(X)B\right)\ .$$
\end{lem}

{\bf Proof} The computation
$$\begin{array}{rl}
\displaystyle \rho(X)(AB)&
\displaystyle =\rho(X)\left(\sum_{\mathbf X,\mathbf Y\in\mathcal X^*}
(A,\mathbf X)(B(\mathbf Y)\mathbf X\mathbf Y\right)\\
&\displaystyle =\rho(X)\left(\sum_{\mathbf X\in\mathcal X^*}
(A,\mathbf X)(B,\mathbf 1)\mathbf X\right)+
A\left(\rho(X)
\left(\sum_{\mathbf X\in\mathcal X^*\setminus 1}(B,\mathbf X)\mathbf X
\right)\right)\\
&\displaystyle =(\rho(X)A)\epsilon(B)+A\left(\rho(X)
\left(\sum_{\mathbf X\in\mathcal X^*}(B,\mathbf X)\mathbf X
\right)\right)\\
&\displaystyle =\epsilon(B)\rho(X)A+A(\rho(X)B)\end{array}$$
shows assertion (i).

The computation
$$0=\rho(X)\mathbf 1=\rho(X)(A^{-1} A)=\epsilon(A)\rho(X)A^{-1}
+A^{-1}(\rho(X) A)$$
shows the identity
$$\rho(X)A^{-1}=-\frac{1}{\epsilon(A)}A^{-1}(\rho(X)A)$$
which implies
$$\begin{array}{l}
\displaystyle 
\rho(X)(A^{-1} B)=\epsilon(B)\rho(X)A^{-1}+A^{-1}(\rho(X)B)
\\
\displaystyle \qquad =
-\frac{\epsilon(B)}{\epsilon(A)}A^{-1}\rho(X)A+A^{-1}\rho(X)B
\end{array}$$
and this proves assertion (ii).
\hfill$\Box$

\begin{prop} \label{propprodclosed}
We have the inclusion
$$\overline{(AB)}\subset A\overline B+\overline A$$
for all $A,B\in\mathbb K\dlangle \mathcal X\drangle$.
\end{prop}

\begin{rem}
The inclusion of Proposition \ref{propprodclosed}  can be strict as shown 
by the example 
$\overline{\left((1-X)\left(\sum_{n=0}^\infty X^n\right)\right)}=
\overline{1}=\mathbb K\not= \mathbb K+\mathbb K X=
(1-X)\ \overline{\frac{1}{1-X}}+\overline{1-X}$.
\end{rem}

{\bf Proof of Proposition \ref{propprodclosed}}
The inclusion $\overline A\subset A\overline B
+\overline A$, together with the formula 
$$\rho(X)(AC)=A(\rho(X)C)+\epsilon(C)\rho(X)A\in A\overline C
+\overline A$$
given by assertion (i) of Lemma \ref{lemrhoXAB} applied to $C\in
\overline B$ shows that
$A\overline B+\overline A$ is recursively closed. Since
$AB\in A\overline B\subset A\overline B+\overline A$, the recursively
closed vector space $A\overline B+\overline A$ contains the recursive
closure $\overline{(AB)}$ of the product $AB$.\hfill$\Box$

For the followgin result, see also Exercise 6, Page 31 of \cite{SaSo}:

\begin{cor} \label{corratclosed}
(i) For all $A,B\in\mathbb K\dlangle\mathcal X\drangle$, we have
$$\dim(\overline{AB})\leq \dim(\overline A)+
\dim(\overline B)\ .$$
for the recursive closure $\overline{AB}$ of the product $AB$.

\ \ (ii) Elements of finite complexity form a subalgebra 
of $\mathbb K\dlangle \mathcal X\drangle$. 
\end{cor}

\begin{prop} \label{propratinverse}
For $A\in\mathbb K\dlangle \mathcal X\drangle^*$ we have the equality
$$\overline{A^{-1}}+\mathbb K= A^{-1}(\overline A+\mathbb K)\ .$$
\end{prop}

\begin{cor} \label{corinvcompl}
For $A\in\mathbb K\dlangle \mathcal X\drangle^*$ we have the equality
$$\dim(\overline{A^{-1}}+\mathbb K)=
\dim(\overline A+\mathbb K)$$
and the inequalities
$$\dim(\overline A)-1\leq 
\dim(\overline{A^{-1}})\leq \dim(\overline A)+1\ .$$
\end{cor}

\begin{rem}
The inequalities of Corollary \ref{corinvcompl} are sharp as shown by 
the example 
$$\dim\left(\overline{\left(\frac{1}{1-X}\right)}\right)=1
\hbox{ and }\dim\left(\overline{1-X}\right)=2\ .$$
\end{rem}

\begin{cor}\label{Schutzenbeasy} Elements of
finite complexity form a rationally closed subalgebra 
of $\mathbb K\dlangle \mathcal X\drangle$ containing the subalgebra
$\mathbb K\dlangle\mathcal X\drangle_{rat}$ of all rational elements.
\end{cor}

{\bf Proof of Theorem \ref{thmSchutz}} Recognisable series
are rational by Proposition \ref{proprecognimpliesrational}.

Rational series are of finite complexity by Corollary
\ref{Schutzenbeasy} and thus recognisable by 
Proposition \ref{proprecog=fincompl}.
\hfill$\Box$
 
{\bf Proof of Proposition \ref{propratinverse}} 
Assertion (ii) of Lemma \ref{lemrhoXAB} implies that $A^{-1}(\overline A
+\mathbb K)$ is recursively closed.
Since it contains $1=A^{-1}A\in A^{-1}\overline A$
and $A^{-1}$ we have the inclusion
$\overline{A^{-1}}+\mathbb K\subset A^{-1}(\overline A+\mathbb K)$.
Exchanging the role of $A$ and $A^{-1}$, we get
$$\overline A+\mathbb K\subset A(\overline{A^{-1}}+\mathbb K)
\subset AA^{-1}(\overline A+\mathbb K)=\overline A+\mathbb K$$
which shows the equality $\overline A+\mathbb K=A\left(\overline{A^{-1}}+
\mathbb K\right)$ equivalent to $\overline{A^{-1}}+\mathbb K=
A^{-1}\left(\overline A+\mathbb K\right)$ after exchange of $A$ and 
$A^{-1}$.\hfill$\Box$

{\bf Proof of Corollary \ref{corinvcompl}} The equality
$\dim(\overline{A^{-1}}+\mathbb K)=\dim(\overline A+\mathbb K)$
follows trivially from Proposition \ref{propratinverse}.

The inequalities follow from the 
inequalities
$$\dim(\overline{A^{-1}})\leq \dim(\overline{A^{-1}}
+\mathbb K)=\dim(\overline A+\mathbb K)\leq \dim(\overline A)+1$$
together with the opposite inequality
$$\dim(\overline{A})\leq \dim(\overline{A^{-1}})+1\ .$$

{\bf Proof of Corollary \ref{Schutzenbeasy}} 
Elements of finite complexity
form clearly a vector space which is a subalgebra of $\mathbb K
\dlangle\mathcal X\drangle$ by assertion (ii) of Corollary 
\ref{corratclosed}. This subalgebra is 
rationally closed by Corollary \ref{corinvcompl}. 

Since this algebra contains obviously the polynomial
subalgebra $\mathbb K\langle \mathcal X\rangle$
of $\mathbb K\dlangle \mathcal X\drangle$, it contains
the rational subalgebra $\mathbb K\dlangle \mathcal X\drangle_{rat}$
of $\mathbb K\dlangle \mathcal X
\drangle$.\hfill$\Box$

\section{Normal forms}\label{secnormal}

This Section introduces normal forms and uses them to solve
enumerative problems over finite fields. Subsections \ref{subsecnormal}
and \ref{subsectreerep} have large overlaps with Chapter II
of \cite{BR} and . Subsections \ref{subsecenum} and 
\ref{sectfastcomp} contain perhaps some new material.

Minimal presentations for
rational series $A\in\mathbb K\dlangle\mathcal X\drangle_{rat}$ 
are not unique but depend on the choice of $A_2,\dots$
completing $A_1=A$ to a basis of $\overline A$. 
A normal form consists of a preferred basis
$A_1=A,A_2,\dots$ of $\overline A$. It selects thus
a unique minimal presentation for every rational series.

This is useful for computational purposes and for solving 
some enumerative problems. 

\subsection{Normal forms and the associated minimal presentations}
\label{subsecnormal}

We suppose henceforth that the finite set $\mathcal X$ is totally 
ordered. We extend the total order of
$\mathcal X$ (right-left) lexicographically to a total order 
of $\mathcal X^*$ by
$$\mathbf Z<\mathbf XX\mathbf Z<\mathbf Y
Y\mathbf X$$ 
for all $\mathbf X,\mathbf Y,\mathbf Z\in\mathcal X^*$
and for all $X,Y\in\mathcal X$ such that $X<Y$.

Given an element $A\in\mathbb K\dlangle\mathcal X\drangle$, 
we consider the (perhaps infinite or empty) increasing sequence
$\mathbf X_1<\mathbf X_2<\dots\subset \mathcal X^*$ constructed 
as follows: If $A=0$ then the associated sequence is empty.
Otherwise, we start with $\mathbf X_1=\emptyset$ and define 
$\mathbf X_{k+1}$ recursively as the smallest element of
$$\{\mathbf X\in \mathcal X^*\ 
\vert\ \rho(\mathbf X)A\not\in\oplus_{j=1}^k\mathbb K
\rho(\mathbf X_j)A\}$$
if this set is nonempty. Otherwise, the 
sequence $\mathbf X_1,\dots$ is the finite sequence
$\mathbf X_1,\mathbf X_2,\dots,\mathbf X_k$.

We call the sequence $\mathbf X_1,\mathbf X_2,\dots$ the {\it 
normal sequence} associated to $A$.

We identify $\mathcal X^*=\{X_1,\dots,X_k\}^*$ with (the
vertices of) the infinite $k-$regular tree rooted  at
the emptyset $\emptyset\in\mathcal X^*$. The $k$ children of a vertex
$\mathbf X\in\mathcal X^*$ are given by
$\mathcal X\mathbf X=\{X_1\mathbf X,X_2\mathbf X,\dots,X_k\mathbf X\}$.
A subset $\mathcal S=\{\mathbf X_1,\mathbf X_2,\dots\}\subset \mathcal X^*$ 
is a {\it subtree} if its vertices form a subtree of $\mathcal X^*$ 
rooted at the root-vertex $\emptyset$ of $\mathcal X^*$.

\begin{prop} \label{propnormalbasis}
(i) The normal sequence $\mathbf X_1,\mathbf X_2,\dots$ associated to $A$
defines a sequence $\rho(\mathbf X_1)A,\rho(\mathbf X_2)A,
\dots\subset \overline A$ of linearly independent elements. They form
a basis of $\overline A$ if $A$ is rational.

\ \ (ii) The elements
$\{\mathbf X_1,\mathbf X_2,\dots\}$ of a non-empty normal
sequence $\mathbf X_1,\mathbf X_2,\dots$ form  
a subtree of $\mathcal X^*$.
\end{prop}

\begin{rem} The normal sequence $\mathbf X_1,
\mathbf X_2,\dots$ is always infinite if 
$A\in \mathbb K\dlangle\mathcal X\drangle$ has infinite complexity. 
The span of
$\rho(\mathbf X_1)A,\rho(\mathbf X_2)A,\dots$ is in general a strict
subspace of $\overline A$ if $\dim(\overline A)=\infty$. 
\end{rem}

Assertion (i) of Proposition \ref{propnormalbasis} selects a
preferred {\it normal basis} 
$\rho(\mathbf X_1)A=A,\rho(\mathbf X_2)A,\dots$
of $\overline A$ for $A\in\mathbb K\dlangle\mathcal X\drangle_{rat}$. 
The corresponding minimal 
presentation is the {\it normal presentation} of $A$.

{\bf Proof of Proposition \ref{propnormalbasis}} Assertion (i)
is obvious by construction of the normal sequence 
$\mathbf X_1,\mathbf X_2,\dots$ associated to $A$.

In order to establish assertion (ii) it is enough to prove that
$\tilde{\mathbf X}\in\{\mathbf X_1,\mathbf X_2,\dots\}$ for 
the immediate ancestor $\tilde{\mathbf X}$ of
$\mathbf X_k=X\tilde{\mathbf X}\in\{\mathbf X_2,\mathbf X_3,\dots\}$.
Since $\tilde{\mathbf X}<\mathbf X_k$, we have either
$\tilde{\mathbf X}\in\{\mathbf X_1,\dots,\mathbf X_{k-1}\}$ 
and we are done or we have
$\rho(\tilde{\mathbf X})A\in\oplus_{j=1}^l\mathbb K
\rho(\mathbf X_j)A$ where $l<k$ is the integer defined by the inequalities
$\mathbf X_1<\mathbf X_2<\dots<
\mathbf X_l<\tilde{\mathbf X}<\mathbf X_{l+1}$. We have thus
$X\mathbf X_1<\dots<X\mathbf X_l<X\tilde{\mathbf X}=\mathbf X_k$
which implies
$$\rho(\mathbf X_k)A=\rho(X\tilde{\mathbf X})A\in
\oplus_{j=1}^l\mathbb K\rho(X\mathbf X_j)A\subset
\oplus_{j=1}^{k-1} \rho(\mathbf X_j)A$$ 
in contradiction with the definition of $\mathbf X_k$.
\hfill$\Box$

\begin{rem} There are many other total orders on $\mathcal X^*$
giving rise to normal forms with good properties: One can consider 
any total order
on $\mathcal X^*$ satisfying $\mathbf X<X\mathbf X$ and
$X\mathbf X<X\mathbf Y$ for all $X\in\mathcal X$ and for all
$\mathbf X,\mathbf Y\in\mathcal X^*$ such that $\mathbf X<\mathbf Y$.
The properties of the lexicographical order will however be needed in
Section \ref{sectfastcomp}.
\end{rem}

\subsection{Tree presentations}\label{subsectreerep}

A subtree $T\subset\mathcal X^*$ of the infinite rooted $k-$regular tree
$\mathcal X^*=\{X_1,\dots,X_k\}^*$ is {\it full} 
if every vertex of $T$ is either a leaf of $T$ or all its $k$ children
are also vertices of $T$. A vertex of the second kind is called an
interior vertex. We denote by $\partial V(T)$ the set of leaves of $T$ and
by $V^\circ (T)$ the set of interior vertices of $T$. 

Let $\mathcal{FFT}(\mathcal X)$ be the set of all finite full
subtrees of $\mathcal X^*$.

{\bf Definition} 
A {\it tree presentation} is a triplet $(T,\epsilon,\mu)$
where $T\subset \mathcal{FFT}(\mathcal X)$ is a finite full subtree of
$\mathcal X^*$ endowed with two maps $\epsilon:V^\circ(T)
\longrightarrow\mathbb K$ and
$\mu:\partial V(T)\times V^\circ(T) \longrightarrow \mathbb K$
such that $\mu(\mathbf X,\mathbf Y)=0$ if $\mathbf X<\mathbf Y$
for $(\mathbf X,\mathbf Y)\in\partial V(T)\times V^\circ(T)$.

A tree presentation encodes a rational series $A=A_\emptyset\in
\mathbb K\dlangle\mathcal X\drangle_{rat}$: Consider the 
recursive presentation with solution given by the series
$A_{\mathbf X},\ \mathbf X\in V^\circ(T)$ indexed by interior vertices
and defined by the equations
$$A_{\mathbf X}=\epsilon(\mathbf X)+\sum_{X\in\mathcal X}
(\rho(X)A_{\mathbf X})\ X,\ \mathbf X\in V^\circ(T)$$
where $\rho(X)A_{\mathbf X}=A_{\mathbf Y}$ if $\mathbf Y=X\mathbf X\in
V^\circ(T)$ and 
$$\rho(X)A_{\mathbf X}=\sum_{\mathbf Y\in V^\circ(T)}
\mu(X\mathbf X,\mathbf Y)\  A_{\mathbf Y}$$
otherwise.

\begin{rem} We have $\rho(\mathbf X)A=A_{\mathbf X}$ for all
$\mathbf X\in V^\circ(T)$ if $A=A_\emptyset$ is defined by 
a tree presentation $(T,\epsilon,\mu)$.
\end{rem}

We call a tree presentation $(T,\epsilon,\mu)$ {\it minimal} if the 
corresponding presentation of $A=A_{\emptyset}\in\mathbb K\dlangle
\mathcal X\drangle_{rat}$ is minimal. 

\begin{prop} \label{propmintreepres}
(i) If $(T,\epsilon,\mu)$ is a tree presentation of 
$A\in\mathbb K\dlangle\mathcal X\drangle_{rat}$ then $V^\circ(T)$
contains all elements $\mathbf X_1,\mathbf X_2,\dots$ of the normal
form of $A$.

\ \ (ii) Every rational series $A\in\mathbb K\dlangle\mathcal 
X\drangle_{rat}$ has a unique minimal tree presentation 
$(T_A,\epsilon,\mu)$. The associated presentation is
the normal presentation of $A$.
\end{prop} 

We call the tree $T_A\in\mathcal{FFT}(\mathcal X)$ 
underlying the normal presentation of $A$ the
{\it minimal tree} of $A$. It has $a=\dim(\overline A)$ interior vertices
given by the elements $\mathbf X_1,\mathbf X_2,\dots$ of the 
normal form associated to $A$ and $1+a(\sharp(\mathcal X)-1)$
leaves given by $\{\emptyset,\mathcal X\mathbf X_1,\dots,\mathcal X
\mathbf X_a\}\setminus\{\mathbf X_1,\dots,\mathbf X_a\}$.

{\bf Proof of Proposition \ref{propmintreepres}}
Suppose that $\mathbf X_j$ is the smallest element of the
normal form $\mathbf X_1,\mathbf X_2,\dots$ associated to $A
\in\mathbb K\dlangle\mathcal X\drangle_{rat}$ which is not an interior 
vertex of the tree $T$ underlying a tree presentation $(T,\epsilon,\mu)$
of $A$. Assertion (ii) of Proposition \ref{propnormalbasis} implies
that $\mathbf X_j$ is of the form $\mathbf X_j=X\mathbf X_{j'}$ for some
$X\in\mathcal X$ and for some element $\mathbf X_{j'}\in
\{\mathbf X_1,\mathbf X_2,\dots,\mathbf X_{j-1}\}\subset V^\circ(T)$.
This shows $\mathbf X_j\in\partial V(T)$ and we have thus 
$$\rho(\mathbf X_j)A\in\sum_{\mathbf Y\in V^\circ(T),\ \mathbf Y<
\mathbf X_j}\mathbb K \rho(\mathbf Y)A=\sum_{i=1}^{j-1}
\mathbb K\rho(\mathbf X_i)A$$
contradicting the inclusion $\mathbf X_j\in\{\mathbf X_1,\mathbf X_2,
\dots\}$. Assertion (i) follows. 

Let $A\in\mathbb K\dlangle\mathcal X\drangle_{rat}$ be a rational series
of complexity $a=\dim(\overline A)$.
Assertion (ii) of Proposition \ref{propnormalbasis}
shows that we can consider the finite full tree 
$T_A\in \mathcal{FFT}(\mathcal X)$ 
with vertices the $1+a\sharp(\mathcal X)$ elements of the form
$$\{\emptyset\}\cup
\{\mathbf X_1,\dots,\mathbf X_a\}\cup\{\mathcal X\mathbf X_1,\dots,
\mathcal X\mathbf X_a\}\subset \mathcal X^*\ .$$
The tree $T_A$ has $a$ interior vertices 
$V^\circ(T_A)=\{\mathbf X_1,\dots,\mathbf X_a\}$ and 
$1+a(\sharp(\mathcal X)-1)$
leaves $\partial V(T_A)\subset \{\emptyset,\mathcal X
\mathbf X_1,\dots,\mathcal X\mathbf X_a\}$.

Interior vertices of $T_A$ are in bijection with the
normal basis $\rho(\mathbf X_1)A,\dots,\rho(\mathbf X_a)A$ of 
$\overline A$ and can thus be endowed by the map 
$\tilde\epsilon(\mathbf X_i)=\epsilon(\rho(\mathbf X_i)A)$
where $\epsilon:\overline A\longrightarrow \mathbb K$ is the 
augmentation map.
For each leaf $\mathbf L=X\mathbf X_i\in\partial V(T_A)$ we consider 
the map
$V^\circ(T_A)\ni\mathbf Y\longmapsto \mu(\mathbf L,\mathbf Y)\in\mathbb K$
defined by the equality 
$$\rho(\mathbf L)A=\sum_{\mathbf Y\in V^\circ(T_A)} 
\mu(\mathbf L,\mathbf Y)\rho(\mathbf Y)A\ .$$
This map satisfies $\mu(\mathbf L,\mathbf Y)=0$ if 
$\mathbf L<\mathbf Y$ by definition of a normal form.
The triplet $(T_A,\tilde\epsilon,\mu)$ is thus a tree presentation 
of some rational element $\tilde A$. 
The associated recursive presentation coincides
by construction with the normal presentation of $A$.
This shows $\tilde A=A$. Minimality of $(T_A,\tilde 
\epsilon,\mu)$ is obvious and unicity follows from assertion (i)
above.
\hfill$\Box$

\subsection{Enumerating elements of given complexity 
in $\mathbb F_q\dlangle X_1,\dots,
X_k\drangle_{rat}$}\label{subsecenum}

The main ingredients for enumerating elements of $\mathbb F_q
\dlangle X_1,\dots,X_k\drangle_{rat}$ according to their complexity
are Proposition \ref{propmintreepres} and the following result:

\begin{prop} \label{propmintreerep}
Let $A\in\mathbb K\dlangle\mathcal X\drangle_{rat}$ be a rational series. 
If $T\in\mathcal{FFT}(\mathcal X)$
contains the minimal tree $T_A$ of $A$, then the set of tree presentations
of $A$ with underlying tree $T$ has a structure of an affine
$\mathbb K-$vectorspace of dimension 
$$\sum_{\mathbf L\in\partial V(T)}\sharp\{
\mathbf Y\in V^\circ(T)\setminus V^\circ(T_A)\ \vert\ \mathbf Y<\mathbf L\}
\ .$$
\end{prop}

{\bf Proof} Let $T\in\mathcal{FFT}(\mathcal X)$ be a tree containing 
the minimal tree $T_A$ of a rational series 
$A\in\mathbb K\dlangle\mathcal X\drangle_{rat}$ with
normal form $\mathbf X_1,\mathbf X_2,\dots,\mathbf X_a$.
Consider the triplet $(T,\tilde \epsilon,\mu)$ with
$\tilde \epsilon:V^\circ(T)\longrightarrow \mathbb K$ defined by
$\tilde \epsilon(\mathbf X)=\epsilon(\rho(\mathbf X)A)$
for every interior vertex $\mathbf X\in V^\circ(T)$ and 
$\mu:\partial V(T)\times V^\circ(T)\longrightarrow \mathbb K$
defined by 
$$\rho(\mathbf L)A=\sum_{i=1}^a \mu(\mathbf L,\mathbf X_i)\rho(\mathbf X_i)A
\in\oplus_{i=1}^a\mathbb K\rho(\mathbf X_i)A=\overline A$$
on $\partial V(T)\times V^\circ(T_A)$ and extended to
$\partial V(T)\times V^\circ(T)$ by setting $\mu(\mathbf L,\mathbf Y)=0$
if $\mathbf Y\in V^\circ(T)\setminus V^\circ(T_A)$. The inclusion
$V^\circ(T_A)\subset V^\circ(T)$ and properties of normal forms show easily
that $(T,\tilde\epsilon,\mu)$ is a tree presentation of $A$.

Since the map $\epsilon:\partial V(T)\longrightarrow \mathbb K$ 
of a tree presentation $(T,\epsilon,\mu)$ of $A$ depends only on $A$ and
the tree $T$ (containing the minimal tree $T_A$ of $A$), an 
arbitrary tree presentation with underlying tree $T$ is of the form
$(T,\tilde \epsilon,\mu')$ for a suitable map $\mu'$ which differs
from $\mu$ by relations among the series $\rho(\mathbf Y)A,\ \mathbf Y
\in V^\circ(T)$. More precisely, the restrictions
$V^\circ(T)\ni\mathbf Y\longmapsto (\mathbf L,\mathbf Y)$
are well-defined up to linear relations in
$$\mathcal S_{<\mathbf L}=\{\rho(\mathbf Y)A\ \vert\ \mathbf Y\in V^\circ(T)
\hbox{ and }\mathbf Y<\mathbf L\}\ .$$
A basis of the vector space spanned by $\mathcal S_{<\mathbf L}$ is given
by
$$\{\rho(\mathbf Y)A\ \vert\ \mathbf Y\in V^\circ(T_A)\hbox{ and }
\mathbf Y<\mathbf L\}\ .$$
The vector space of linear relations among elements of $\mathcal S_{<\mathbf L}
$ is thus of dimension 
$$\{\mathbf Y\in V^\circ(T)\setminus V^\circ(T_A)\ \vert \ 
\mathbf Y<\mathbf L\}\ .$$
A summation over all leaves $\mathbf L\in\partial V(T)$ shows the result.
\hfill$\Box$

We endow the set $\mathcal{FFT}(\mathcal X)$ of finite full 
subtrees in $\mathcal X^*$ with the partial order given by inclusion:
$T'<T$ for $T',T\in\mathcal{FFT}(\mathcal X)$ if 
$V^\circ(T')$ is a strict subset of $V^\circ(T)$.

For $T\in\mathcal{FFT}(\mathcal X)$, we define $E_{T}(q)\in
\mathbb Z[q]$ recursively by 
$$E_{T}(q)=q^{\sharp V^\circ(\mathbf T)}\prod_{\mathbf L\in\partial V(T)}
q^{\sharp\{\mathbf X\in V^\circ(T)\ \vert\ \mathbf X<\mathbf L\}}
-C_{T}(q)$$
where 
$$C_{T}(q)=\sum_{T'\in\mathcal{FFT}(\mathcal X),
\ T'<T}E_{T'}(q)\prod_{\mathbf L\in\partial
V(T)}q^{\sharp\{\mathbf X<\mathbf L\ \vert\ 
\mathbf X\in V^\circ(T)\setminus V^\circ(T')\}}\ .$$

\begin{prop} \label{propenumlin}
For $\mathbb K=\mathbb F_q$ a finite field with $q$
elements, the integer $E_{\mathbf T}(q)$ is the number of elements
in $\mathbb F_q\dlangle\mathcal X\drangle$ with minimal tree $T$.
\end{prop}

We consider a sequence $E_0(q),E_1(q),\dots\subset \mathbb Z[q]$ 
defined by 
$$E_n(q)=\sum_{\mathbf T\in\mathcal{FFT}(\mathcal X),\ 
\sharp V^\circ(\mathbf T)=n}
E_{\mathbf T}(q)\ .$$

\begin{cor} (i) The number of rational series with complexity $n$ 
in $\mathbb F_q\dlangle\mathcal X\drangle_{rat}$ is given by
$E_n(q)$.

\ \ (ii) We have 
$$\sum_{A\in\mathbb F_q\dlangle\mathcal X\drangle_{rat}}
t^{\dim(\overline A)}=\sum_{n=0}^\infty E_n(q)t^n\ .$$
\end{cor}

{\bf Proof of Proposition \ref{propenumlin}} Given a finite full tree 
$T\in\mathcal{FFT}(\mathcal X)$ there are $q^{\sharp V^\circ(T)}$
possible choices for the map $\epsilon:V^\circ(T)\longrightarrow
\mathbb F_q$ and 
$$\prod_{\mathbf L\in \partial V(T)}q^{\sharp\{\mathbf X\in
V^\circ(T)\ \vert\ \mathbf X<\mathbf L\}}$$
possible choices for the map $\mu:\partial V(T)\times
V^\circ(T)\longrightarrow \mathbb K$ 
giving rise to a tree presentation $(T,\epsilon,\mu)$
with underlying tree $T$. 

Proposition \ref{propmintreerep} shows that
a rational series $A$ with minimal tree $T_A$ contained in $T$ is 
has exactly 
$$\prod_{\mathbf L\in\partial V(T)}q^{\sharp\{\mathbf X\in V^\circ(T)
\setminus V^\circ(T')\ \vert \ \mathbf X<\mathbf L\}}$$
tree presentations with underlying tree $T$.
This leads to the correction $C_T(q)$ and shows the result.
\hfill$\Box$

\begin{exple} One gets easily $E_n(q)=q^{2n}-q^{2n-1}$ in the commutative case
$\mathcal X=\{X\}$ involving a unique variable.
\end{exple}

\begin{exple} Describing an element $T\in\mathcal{FFT}(\mathcal X)$
by the set of its vertices, the first polynomials $E_T(q)$ for 
$\mathcal X=\{X,Y\}$ are: 
$$\begin{array}{l}
\displaystyle
E_{\{\emptyset\}}=1\\
\displaystyle
E_{\{\emptyset,X,Y\}}=q^3-q^2\\
\displaystyle 
E_{\{\emptyset,X,X^2,YX,Y\}}=q^8-2q^6+q^5\\
\displaystyle 
E_{\{\emptyset,X,Y,YX,Y^2\}}=q^7-2q^5+q^4\\
\displaystyle
E_{\{\emptyset,X,X^2,X^3,YX^2,YX,Y\}}=q^{15}-2q^{12}-q^{11}+3q^{10}-q^9\\
\displaystyle
E_{\{\emptyset,X,X^2,YX,XYX,Y^2X,Y\}}=q^{14}-2q^{11}-q^{10}+3q^9-q^8\\
\displaystyle
E_{\{\emptyset,X,X^2,YX,Y,XY,Y^2\}}=q^{13}-q^{11}-2q^{10}+q^9+2q^8-q^7\\
\displaystyle
E_{\{\emptyset,X,Y,XY,X^2Y,YXY,Y^2\}}=q^{13}-2q^{10}-q^9+3q^8-q^7\\
\displaystyle
E_{\{\emptyset,X,Y,XY,Y^2,XY^2,Y^3\}}=q^{12}-2q^9-q^8+3q^7-q^6
\end{array}$$
The algebra $\mathbb F_q\dlangle X,Y\drangle_{rat}$ contains thus 
$E_1(q)=q^3-q^2$ elements of complexity $1$, 
$$E_2(q)=q^8+q^7-2\,q^6-q^5+q^4$$
elements of complexity $2$ and 
$$E_3(q)=q^{15}+q^{14}+2\, q^{13}-q^{12}-4\,q^{11}-2\,q^{10}+3\,q^8+q^7-q^6$$
elements of complexity $3$.

The computation of 
$$\begin{array}{l}
\displaystyle
E_4(q)={q}^{24}+{q}^{23}+2\,{q}^{22}+3\,{q}^{21}+{q}^{20}-7\,{q}^{18}-8\,{q}^
{17}-3\,{q}^{16}+\\
\displaystyle\qquad\qquad +6\,{q}^{14}+5\,{q}^{13}
+3\,{q}^{12}-4\,{q}^{10}-{q}^{
9}+{q}^{8}\end{array}$$
is already tedious and motivates the approach given below.
\end{exple}

\subsection{Fast computation of $E_n(q)$ for $\mathcal X=\{X,Y\}$}
\label{sectfastcomp}

This section contains formulae (without proofs)
for efficient computations of the
polynomials $E_n(q)$ in the case where $\mathcal X=\{X,Y\}$
has two elements. 

We set $w_0=1$ and define $w_{n+1}(q)\in\mathbb Z[q]$ recursively by the
formula
$$w_{n+1}=q^{3+n}\sum_{j=0}^nq^{j(n+1-j)}w_jw_{n-j}\ .$$

\begin{prop} We have
$$w_n(q)=q^n\sum_{T\in\mathcal{FFT}(\mathcal X),\ V^\circ(T)=n}
\prod_{\mathbf L\in \partial V(T)}q^{\sharp\{
\mathbf X\in V^\circ(T)\ \vert\ \mathbf X<\mathbf L\}}\ .$$
\end{prop}

The proof is left to the reader.

Working with the lexicographic order on $\mathcal X^*$,
a little thought shows that the contribution of a tree 
$T'\in\mathcal{FFT}(\mathcal X)$ with $j=\sharp(V^\circ(T'))$ 
interior vertices to 
$$\sum_{T\in\mathcal{FFT}(\mathcal X),\ V^\circ(T)=n}
C_T(q)$$
is of the form $p_{j,n}(q)E_{T'}(q)$ with 
$p_{j,n}\in\mathbb Z[q]$ a polynomial depending only on
$j=\sharp(V^\circ(T'))$ and $n$.

This implies the recursive formula
$$E_n(q)=w_n(q)-\sum_{j=0}^{n-1}p_{j,n}(q)E_j(q)$$
for the polynomials $E_0(q),E_1(q),\dots$.

A little work shows that the polynomials $p_{j,n}(q)$ are given
by $p_{0,n}(q)=q^{-n}w_n(q)$ for $j=0$. The remaining values
$p_{j,n}(q),\ 1\leq j\leq n$ can be recursively computed using
the formulae
$$\begin{array}{rcl}
\displaystyle p_{j,n}(q)&
\displaystyle =&
\displaystyle
\sum_{h=0}^{n-j}q^{(h+1)(n-j-h)}q^{-h}w_h(q)p_{j-1,n-1-h}(q)\\
&\displaystyle =&
\displaystyle
\sum_{h=0}^{n-j}q^{h(n-1-h)}w_h(q)p_{j-1,n-1-h}(q)\end{array}$$
leading to the same result.

This shows that the computation of $E_n(q)$ can be done in polynomial
time (with respect to $n$).

For $q=2$, the first few coefficients of the 
series $\sum_{n=0}^\infty E_n(2)t^n$ are
$$\begin{array}{l}
1+4\,t+240\,{t}^{2}+52032\,{t}^{3}+37961472\,{t}^{4}+95557604352\,{t}^
{5}\\
+873176389545984\,{t}^{6}+30234012628981334016\,{t}^{7}\\+
4073184753921806027390976\,{t}^{8}+2164965110784257951109280432128\,{t
}^{9}\\
+4571419424684923104187906920444592128\,{t}^{10}\\+
38479163698041617829387740718124411857666048\,{t}^{11}\\+
1293355066072995022042530447708918263083390363238400\,{t}^{12}\\+
173739578583285839772280634310511087695154611244324192518144\,{t}^{13}
\end{array}$$

\begin{rem} Similar more complicated
formulae for polynomial time algorithms 
exist for arbitrary finite sets $\mathcal X$.
\end{rem}

\begin{rem} Similar (but slightly trickier) arguments give the
number $\tilde E_n(q)$ of polynomials of complexity $n$ in
$\mathbb F_q\langle X,Y\rangle$ by the formula
$$\tilde E_n(q)=\tilde w_n(q)-\sum_{j=0}^{n-1}\tilde p_{j,n}(q)\tilde 
E_j(q)$$
where 
$$\tilde w_0(q)=1,\ \tilde
w_{n+1}=q\sum_{j=0}^nq^{j(n+1-j)}\tilde w_j\tilde w_{n-j},\ n\geq 0$$
and where
$\tilde p_{0,n}(q)=q^{-n}\tilde w_n(q)$,
$$\begin{array}{rcl}
\displaystyle \tilde p_{j,n}(q)&
\displaystyle =&
\displaystyle
\sum_{h=0}^{n-j}q^{(h+1)(n-j-h)}q^{-h}\tilde w_h(q)\tilde p_{j-1,n-1-h}(q)\\
&\displaystyle =&
\displaystyle
\sum_{h=0}^{n-j}q^{h(n-1-h)}\tilde w_h(q)\tilde p_{j-1,n-1-h}(q)\end{array}$$
for $1\leq j\leq n$.

The generating series $\sum_{n=0}^\infty E_n(2)t^n$ starts as
$$\begin{array}{l}
1+t+6\,{t}^{2}+72\,{t}^{3}+1776\,{t}^{4}+89280\,{t}^{5}+9065472\,{t}^{
6}+1850148864\,{t}^{7}\\
+757046525952\,{t}^{8}+620298979246080\,{t}^{9}+
1017126921430892544\,{t}^{10}\\
+3336658943759213395968\,{t}^{11}+
21894988380633154342354944\,{t}^{12}\\
+287369531352172835754234347520\,{
t}^{13}\\
+7543680108676972971562235527692288\,{t}^{14}\\+
396062820851396884301553848136757149696\,{t}^{15}\\+
41589051965658313888146456051766022098649088\,{t}^{16}\\+
8734258246436387382993841213619134491337634611200\,{t}^{17}\\+
3668631292951234310193522386177325845447083530708320256\,{t}^{18}
\end{array}$$
and these formulae have again generalisations to an arbitrary number of
variables.
\end{rem}

\section{Saturation level}\label{sectsaturlevel}

We denote by $J_n:\mathbb K\dlangle
\mathcal X\drangle\longrightarrow \mathbb K\langle
\mathcal X\rangle$
the $n-$jet, namely the linear projection defined by 
$$\mathbb K\dlangle\mathcal X\drangle\ni
A\longmapsto
J_n(A)=\sum_{\mathbf X\in\mathcal X^{\leq n}} (A,\mathbf X)\mathbf X$$
where the summation is over all words $\mathbf X\in\mathcal X^*
\setminus\mathfrak m^{n+1}$ of length $\leq n$.

Given a subspace 
$\mathcal A\subset \mathbb K\dlangle \mathcal X\drangle$,
we denote by $J_n(\mathcal A)\subset \mathbb K\langle\mathcal X\rangle$
its image under the projection $J_n$ and by 
$$\mathcal K_n(\mathcal A)=
\mathop{ker}(J_n)\cap \mathcal A=\{A\in\mathcal A\ \vert
\ J_n(A)=0\}\subset \mathcal A$$
the kernel of the projection $J_n$ restricted to $\mathcal A$.
The vector spaces $\mathcal K_n(\mathcal A),\mathcal A$ and $J_n(\mathcal A)$
are related by the exact sequence
$$0\longrightarrow\mathcal K_n(\mathcal A)\longrightarrow \mathcal A
\longrightarrow J_n(\mathcal A)\longrightarrow 0\ .$$

A useful tool for computations with rational series
is the {\it saturation degree} or {\it saturation level}:
Given a recursively closed
subspace $\mathcal A\subset \mathbb K\dlangle \mathcal X\drangle$, the
saturation level of $\mathcal A$ is the smallest element $N\in
\mathbb N\cup\{\infty\}$ 
such that $\mathcal K_N(\mathcal A)=\mathcal K_{N+1}
(\mathcal A)$.

\begin{prop} If a recursively closed subspace
$\mathcal A=\mathbb K\dlangle \mathcal X\drangle$
has saturation degree $N$ then 
$J_N: \mathcal A\longrightarrow \pi_N(\mathcal A)\subset
\mathbb K\langle \mathcal X\rangle$
is one-to-one.
\end{prop}

{\bf Proof}
We have by definition the inclusions
$$\mathcal K_0(\mathcal A)=\mathcal A\cap\mathfrak m
\supset \mathcal K_1(\mathcal A)\supset\mathcal K_2
(\mathcal A)\supset \dots \ .$$
The equality $\mathcal K_N(\mathcal A)=\mathcal K_{N+1}(\mathcal A)$,
together with the obvious 
inclusions $\rho(X)\left(\mathcal K_{l+1}(\mathcal A)\right)\subset
\mathcal K_l(\mathcal A)$ for all $l\in\mathbb N$ and for all $X\in
\mathcal X$,
shows that $\mathcal K_N(\mathcal A)\subset \mathcal A$ 
is a recursively closed subspace of $\mathcal A\cap\mathfrak m$. 
This shows $\mathcal K_N(\mathcal A)=\{0\}$
since $\{0\}$ is the only recursively closed subspace of 
$\mathfrak m$. Indeed, if  $B$ is non-zero,
there exists a monomial $\mathbf X\in\mathcal X^*$
whose coefficient $(B,\mathbf X)$ is non-zero.
This implies $\epsilon(\rho(\mathbf X)B)=(B,\mathbf X)\not= 0$.
We have thus $\rho(\mathbf X)B\not\in \mathfrak m$ implying
$\overline B\not\subset \mathfrak m$.
\hfill$\Box$

The saturation level is useful for proving the following result,
cf also Proposition 3.1, Page 46 of \cite{BR}:

\begin{prop} \label{propfinmonoid}
The following statements are equivalent:

\ \ (i) The shift monoid $\rho_{\overline A}(\mathcal X^*)\subset
\hbox{End}(\overline A)$ of $A\in\mathbb K\dlangle \mathcal X\drangle$
is finite.

\ \ (ii) $A\in\mathbb K\dlangle\mathcal X\drangle$ is rational 
and has all its coefficients in a finite subset of $\mathbb K$.
\end{prop}

{\bf Proof} Suppose that $A$ has a finite shift monoid $\rho_{\overline A}
(\mathcal X^*)$. This implies that $\rho_{\overline A}(\mathcal X^*)A$
is finite and $\overline A$ is finite-dimensional. 
The set $\{(A,\mathbf X)\ \vert\ \mathbf X\in\mathcal X^*\}$
of coefficients of $A$ is thus given by the finite 
set $\{\epsilon(\rho(\mathbf X)A)\ \vert
\ \mathbf X\in\mathcal X^*\}$.
This shows that (i) implies (ii).

If $A\in\mathbb K\dlangle\mathcal X\drangle$ is rational,
its recursive closure $\overline A$ has finite saturation level $N$
and we have a faithfull map $J_N:\overline A\longrightarrow
\mathbb K\langle\mathcal X\rangle$ into the finite-dimensional 
vector space of non-commutative polynomials of degree $\leq N$.
If all coefficients of $A$ belong to a finite subset $\mathcal F\subset 
\mathbb K$ of $\mathbb K$, the image $J_N(\rho(\mathcal X^*)A)$
is contained in the finite set of non-commutative polynomials of degree
$\leq N$ with coefficients in $\mathcal F$. Since $J_N$ is a faithful on
$\overline A$, the orbit $\rho(\mathcal X^*)A$ is finite.
The shift monoid $\rho_{\overline A}(\mathcal X)$ of $A$
is thus finite since it has a faithful action on the
finite set $\rho(\mathcal X^*)A$.
This shows that (ii) implies (i).\hfill$\Box$

\subsection{Algorithmical aspects}

The properties of the saturation level imply the existence
of finite algorithms for all operations in the rationally closed
algebra $\mathbb K\dlangle \mathcal X\drangle_{rat}$. Two rational
elements $A,B\in\mathbb K\dlangle \mathcal X\drangle_{rat}$
described by finite presentations can be compared,
added and multiplied using only finitely many arithmetical 
operations in $\mathbb K$. 
Similarly, the computation of $A^{-1}$ for an invertible element 
$A\in\mathbb K\dlangle \mathcal X\drangle_{rat}^*$
uses also only a finite number of operations in $\mathbb K$. 

Indeed, the formulae 
$$\begin{array}{l}
\displaystyle 
\rho(X)(A+B)=\rho(X)A+\rho(X)B,\\
\displaystyle
\rho(X)(A_iB_j)=A_i\rho(X)B_j+\epsilon(B_j)\rho(X)A_i,\\
\displaystyle
\rho(X)A^{-1}=-1/\epsilon(A)\ A^{-1}\left(\rho(X)A\right),\\
\displaystyle
\rho(X)(A^{-1}A_i)=A^{-1}\left(-\epsilon(A_i)/\epsilon(A)\ \rho(X)A
+\rho(X)A_i\right)\end{array}$$
for $A,B\in\mathbb K\dlangle\mathcal X\drangle_{rat}$ (with $A$ invertible
for the last two formulae), enables us easily to write down recursive
presentations for
$A+B,AB$ and $A^{-1}$, given recursive presentations of $A$ and $B$. 
Computing the saturation level allows then to compute associated
minimal (or normal) presentations by removing first linearly dependent 
elements and by computing then the exact image of the recursive 
closure of $A+B,AB$ or $A^{-1}$.

\begin{exple} Consider the series $A=A_1=1/(1-XY)$ of Example
\ref{expleunsurunmoinsxy} defined by the recursive presentation
$$\left\lbrace\begin{array}{l}
\displaystyle A_1=1+A_2Y\\
\displaystyle A_2=A_1X\end{array}\right.$$
Setting $B_1=A^{-1},B_2=A^{-1}A_1,B_3=A^{-1}A_2$, we have
$$\epsilon(B_1)=\epsilon(A^{-1})=1/\epsilon(A)=1,\quad
\epsilon(B_2)=\epsilon(A^{-1}A_1)=1,\quad \epsilon(B_3)=0$$
and
$$\begin{array}{l}
\displaystyle \rho(X)B_1=-A^{-1}\rho(X)A=-A^{-1}0=0,\\
\displaystyle \rho(Y)B_1=-A^{-1}\rho(Y)A=-A^{-1}A_2=-B_3,\\
\displaystyle \rho(X)B_2=A^{-1}(-\rho(X)A+\rho(X)A_1)=0,\\
\displaystyle \rho(Y)B_2=A^{-1}(-\rho(Y)A+\rho(Y)A_1)=0,\\
\displaystyle \rho(X)B_3=A^{-1}(-0\rho(X)A+\rho(X)A_2)=B_2,\\
\displaystyle \rho(Y)B_3=A^{-1}(-0\rho(Y)A+\rho(Y)A_2)=0\end{array}$$
leading to the presentation
$$B_1=1-B_3Y,\quad B_2X=1,\quad B_3=B_2X$$
(which is already minimal) 
and showing $B_1=1-B_3Y=1-(B_2X)Y=1-XY$ as expected.
\end{exple}

\begin{rem} The saturation level, although sometimes useful, is by 
no means absolutely necessary for dealing with computational aspects of 
$\mathbb K\dlangle\mathcal X\drangle_{rat}$. The features described 
in Sections \ref{secreducrec}, \ref{secminimal} (and Section 
\ref{subsecnormal} when dealing with comparisons) can be used as a substitut.
\end{rem}

\section{The metric group $S\mathbb K\dlangle \mathcal X\drangle_{rat}^*$}
\label{sectmetricgroup}

An {\it oriented 
norm} on a group $\Gamma$ with identity $e$ is an
application $\parallel \quad \parallel_o:
\Gamma\setminus \{e\}\longrightarrow \mathbb R^*_+=\{x\in\mathbb R\ 
\vert\ x>0\}$
which satisfies the {\it triangle inequality} 
$\parallel \gamma\delta\parallel_o\leq \parallel\gamma\parallel_o+
\parallel\delta\parallel_o$ for all
$\gamma,\delta\in\Gamma$. We extend the oriented norm
$\parallel\quad\parallel_o$ to $\Gamma$ by
setting $\parallel e\parallel_o=0$.

A {\it norm} on $\Gamma$ is an oriented norm which is
{\it symmetric}: $\parallel \gamma\parallel_o=\parallel\gamma^{-1}\parallel_o$
for all $\gamma\in\Gamma$.
Every oriented norm gives rise to a norm $\parallel
\gamma\parallel=\parallel\gamma\parallel_o+\parallel\gamma^{-1}\parallel_o$.
A norm turns the group $\Gamma$ into a homogeneous
metric space by considering the distance
$$d(\gamma,\delta)=d(\beta\gamma,\beta\delta)=\parallel\gamma^{-1}\delta
\parallel$$
for $\beta,\gamma,\delta\in\Gamma$. 
In the sequel a {\it metric group} $(\Gamma,\parallel\quad\parallel)$
is a group $\Gamma$ endowed with a norm $\parallel\quad \parallel$.

\begin{thm} \label{propmultlength} The application 
$$A\longmapsto \parallel A\parallel=
\dim(\overline A+\mathbb K)-1$$
defines a norm on the special rational 
group $S\mathbb K\dlangle\mathcal X\drangle_{rat}^*$.
\end{thm}

{\bf Proof} Consider $A\in S\mathbb K\dlangle \mathcal X\drangle_{rat}^*$.
The identity $\dim(\overline A+\mathbb K)-1=0$ implies
$A=1$ and shows $\parallel A\parallel \geq 1$
if $A\not=1$.

In order to establish the triangle inequality,
we consider $A,B\in S\mathbb K\dlangle\mathcal X\drangle_{rat}^*$.
The main tool is the inclusion 
$\overline{AB}\subset A\overline B+\overline A$ of Proposition 
\ref{propprodclosed}. 

If none of $\overline A,\overline B$ contains $\mathbb K$ then
$$\begin{array}{l}
\displaystyle \parallel AB\parallel=\dim(\overline{AB}+\mathbb K)-1\leq
\dim(A\overline B+\overline A+\mathbb K)-1\leq\\
\displaystyle \leq \dim(\overline A)+
\dim(\overline B)=
\dim(\overline A+\mathbb K)-1+
\dim(\overline B+\mathbb K)-1\\
\displaystyle =\parallel A\parallel+\parallel B\parallel\ .\end{array}$$

If $\mathbb K\subset \overline A$ and $\mathbb K\not\subset \overline
B$ then $A\overline B+\overline A+\mathbb K=A\overline B+\overline A$
and we have
$$\begin{array}{l}
\displaystyle \parallel AB\parallel=\dim(\overline{AB}+\mathbb K)-1\leq
\dim(A\overline B+\overline A)-1\leq\\
\displaystyle \leq \dim(\overline A)+\dim
(\overline B)-1=\dim(\overline A+\mathbb
K)-1+\dim(\overline B+\mathbb K)-1\\
\displaystyle =\parallel A\parallel+\parallel B\parallel\ .\end{array}$$

If $\mathbb K\not\subset \overline A$ and $\mathbb K\subset \overline
B$ then $A\in A\overline B\cap \overline A$ and $A\overline
B+\overline A$ is of dimension at most $\dim(\overline A)+
\dim(\overline B)-1$. This implies
$$\begin{array}{l}
\displaystyle \parallel AB\parallel=\dim(\overline{AB}+\mathbb K)-1\leq
\dim(A\overline B+\overline A+\mathbb K)-1\leq
\dim(A\overline B+\overline A)\leq\\
\displaystyle \leq \dim(\overline A)+\dim
(\overline B)-1=\dim(\overline A+\mathbb
K)-1+\dim(\overline B+\mathbb K)-1\\
\displaystyle =\parallel A\parallel+\parallel B\parallel .\end{array}$$

If $\mathbb K\subset \overline A\cap \overline B$ then
$A\in A\overline B\cap \overline A$ and the dimension of
$A\overline B+\overline A+\mathbb K=A\overline B+\overline A$
is at most $\dim(\overline A)+
\dim(\overline B)-1$. We have thus
$$\begin{array}{l}
\displaystyle \parallel AB\parallel=\dim(\overline{AB}+\mathbb K)-1\leq
\dim(A\overline B+\overline A)-1\leq\\
\displaystyle \leq \dim(\overline A)+\dim
(\overline B)-2=\dim(\overline A+\mathbb
K)-1+\dim(\overline B+\mathbb K)-1\\
\displaystyle \parallel A\parallel +\parallel B\parallel\end{array}$$
which ends the proof of the triangle inequality.

The identity $\parallel A\parallel=\parallel A^{-1}\parallel$
for $A\in S\mathbb K\dlangle \mathcal X\drangle_{rat}^*$ 
follows from the equality $\dim(\overline A+
\mathbb K)=\dim(\overline{A^{-1}}+\mathbb K)$
of Corollary \ref{corinvcompl}.\hfill$\Box$

\begin{rem} The metric group $S\mathbb K\dlangle 
\mathcal X\drangle_{rat}^*$ described by Theorem 
\ref{propmultlength} is a non-commutative analogue of the abelian metric
group of rational fractions in commuting
variables evaluating to $1$ at the origin with norm given by 
$$\parallel f/g\parallel
=\max(\mathop{deg}(f),\mathop{deg}(g))$$ 
where $f/g$ is reduced expression, see also Example 
\ref{remdimratfrac}. This group is of course the free abelian group
on all irreducible monic polynomials over $\mathbb K$. In the
case of one variable with $\mathbb K$ algebraically closed, 
these generators are
all affine polynomials of the form $1+\lambda x,\lambda\in \mathbb K^*$ 
and the considered norm is simply the word length
with respect to the (infinite) symmetric generating system
$\{(1-\lambda x)/(1-\mu x)\}_{\lambda,\mu\in\mathbb K}$. 
\end{rem}

\subsection{The Magnus representation of the free group}

The {\it Magnus representation} 
is the representation of the free group $F_k=\langle
g_1,\dots,g_k\rangle$ on $k$ generators defined by $\mu(g_j)=1+X_j
\in S\mathbb K\dlangle X_1,\dots,X_k\drangle_{rat}^*$ 
(see for instance Th\'eor\`eme 1
of Chapitre II, \S 5 in \cite{Bour23}).  

Recall that every element $g$ of the free group $F_k=\langle g_1,\dots,
g_k\rangle$ has a unique reduced
expression $g=g_{i_1}^{\alpha_1}\cdots g_{i_m}^{\alpha_m}$
with indices $i_j\not=i_{j+1}$ in $\{1,\dots,k\}$
and exponents $\alpha_1,\dots,\alpha_m\in\mathbb Z\setminus\{0\}$.
The function $g\longmapsto
\parallel g\parallel=\sum_{j=1}^m \vert \alpha_j\vert$ defined by the
{\it length} $\vert \alpha_1\vert +\dots+\vert \alpha_m\vert$
of the reduced expression for $g\in F_k$ defines a length function 
on $F_k$. This length function coincides with the combinatorial 
length function on the Cayley graph (given by the
infinite $2k-$regular tree) of $F_k$ with respect to the free
symmetric generating set $\{g_1^{\pm 1},\dots,g_k^{\pm 1}\}$.

\begin{thm} \label{thmmagnusmetr}
Let $g=g_{i_1}^{\alpha_1}\cdots g_{i_m}^{\alpha_m}\in F_k$
be a reduced word. Then
$$\parallel \mu(g)\parallel=-c+\sum_{j=1}^m \vert \alpha_j\vert$$
where 
$$c=\sharp\{1\leq j<m\ \vert\ \alpha_j>0,\ \alpha_{j+1}<0\}\ .$$

Otherwise stated, the norm on $\mu(F_k)$
is the norm on $F_k$ with respect to the symmetric 
generating system $g_1^{\pm 1},\dots,g_k^{\pm 1},
g_ig_j^{-1}$ for $i\not= j$,\ $1\leq 
i,j\leq k$.
\end{thm}

\begin{rem} In particular, Theorem \ref{thmmagnusmetr} shows that the 
length of $\mu(g)\in S\mathbb K\dlangle\mathcal X\drangle_{rat}^*$
is independent of $\mathbb K$.
\end{rem}

The following result is well-known, see for example Th\'eor\`eme 1,
of Page 46 in \cite{Bour23} for a more general statement:

\begin{cor} \label{cormagnus}
The Magnus representation is faithful.
\end{cor}

\begin{rem} Theorem \ref{thmmagnusmetr} implies that 
the Cayley graph of the free monoid generated by $1+X_1,\dots,
1+X_k$ is a rooted $k-$regular tree which embedds isometrically into
the metric group $S\mathbb K\dlangle \mathcal X\drangle_{rat}^*$.
See Section \ref{subsecmetrmon} for a generalisation.
\end{rem}

{\bf Proof of Corollary \ref{cormagnus}} If $g=g_{i_1}^{\alpha_1}
\cdots g_{i_m}^{\alpha_m}$ is a
reduced non-trivial word of the free group $F_k$, then 
$$c=\sharp\{1\leq j<m\ \vert\ \alpha_j>0,\ \alpha_{j+1}<0\}\leq
\frac{m}{2}\leq 
\frac{1}{2}\sum_{j=1}^m \vert \alpha_j\vert<\sum_{j=1}^m \vert 
\alpha_j\vert\ .$$
Theorem \ref{thmmagnusmetr} shows thus
$\parallel \mu(g)\parallel \geq 1$ which implies $\mu(g)\not=1$.
\hfill $\Box$

{\bf Proof of Theorem \ref{thmmagnusmetr}} The easy computations
$$\begin{array}{l}
\displaystyle \rho(X)\left((1+X)\frac{1}{1+Y}\right)=1\\
\displaystyle \rho(Y)\left((1+X)\frac{1}{1+Y}\right)=-(1+X)\frac{1}{1+Y}
\end{array}$$
show that the series $(1+X)\frac{1}{1+Y},\ X,Y$ two distinct 
elements of $\mathcal X$,
have norm $1$. The triangle inequality implies thus
$$\parallel \mu(g)\parallel\leq -c+\sum_{j=1}^m \vert \alpha_j\vert\ .$$

In order to prove the opposite inequality $\parallel
\mu(g)\parallel \geq -c+\sum_{j=1}^m \vert\alpha_j\vert$, 
we rewrite $g$ as a word $g=w_1w_2\cdots w_l$ of length 
$l=-c+\sum_{j=1}^m \vert\alpha_j\vert$ with respect to 
the symmetric generating set
$\mathcal S_k=\{g_i^{\pm 1},g_ig_j^{-1}\}_{1\leq i\not= j\leq n}$.

Setting $\mathcal V_0=\mathbb K$ and $$\mathcal V_{s}=
\mathcal V_{s-1}+\mathbb K\mu(w_1\cdots w_s)\subset
\mathbb K\dlangle\mathcal X\drangle_{rat}^*$$
for $s=1,\dots,l$, we have the following result:

\begin{lem} \label{lemmagnus}
We have $\mathcal V_s=\overline{\mu(w_1\cdots w_s)}+\mathbb K$
and $\dim(\mathcal V_s)=s+1$ for all $s\in\{0,\dots,l\}$.
\end{lem}

This shows $\parallel \mu(g)\parallel=
\dim(\mathcal V_l)-1=l$ and ends the proof of Theorem \ref{thmmagnusmetr}.
\hfill$\Box$

{\bf Proof of Lemma \ref{lemmagnus}} 
Writing $W_j=\mu(w_j)$ we remark that $\rho(X)W\in\{0,1,W\}$
for $W\in \mu(\mathcal S_k)$. Lemma \ref{lemrhoXAB} 
and induction on $s$ imply that
all vector spaces $\mathcal V_0,\dots,\mathcal V_l$ are recursively
closed.

Lemma \ref{lemmagnus} holds clearly for  $s=0$ and $s=1$. 
We prove it by induction on $s$: Consider $W_1\dots W_{s+1}$
for $s\geq 1$. If $W_{s+1}=(1+X)$ or $W_{s+1}=(1+X)/(1+Y)$ with 
$X\not= Y$, then
$$\rho(X)(W_1\cdots W_{s+1})\equiv W_1\cdots W_s\pmod{\mathcal V_{s-1}}$$
and 
$$\rho(X)\mathcal V_s\subset\mathcal V_{s-1}$$
since $W_s\not\in\{1/(1+X),(1+Z)/(1+X)\}$.

If $W_{s+1}=1/(1+X)$ then $W_s\in \{1/(1+X),1/(1+Y),
(1+Y)/(1+X),(1+Z)/(1+Y)\}$ with $X\not=Y$ and $Y\not=Z$. 

If 
$W_s\in \{1/(1+X),(1+Y)/(1+X)\}$ we have
$$(\rho(X)-1)(W_1\cdots W_{s+1})\equiv W_1\cdots W_s\pmod{
\mathcal V_{s-1}}$$
and 
$$(\rho(X)-1)\mathcal V_s\subset \mathcal V_{s-1}\ .$$

If $W_s\in\{1/(1+Y),(1+Z)/(1+Y)\}$
we have
$$(\rho(X)+\rho(Y)-1)
(W_1\cdots W_{s+1})=W_1\cdots W_s\pmod{\mathcal V_{s-1}}$$
and
$$(\rho(X)+\rho(Y)-1)\mathcal V_s\subset \mathcal V_{s-1}\ .$$
There exists thus always an element $R_s\in\mathbb K[\rho(\mathcal X)]$
such that $R_s(W_1\cdots W_{s+1})\equiv W_1\cdots W_s\pmod{\mathcal
V_{s-1}}$ and $R_s\mathcal V_s\subset \mathcal V_{s-1}$.
Setting $\tilde W=W_1\dots W_{s+1}$, the induction hypothesis
shows that the $(s+2)$ elements 
$$\tilde W,R_s\tilde W,R_{s-1}R_s\tilde W,\dots,R_1\cdots
R_s\tilde W,1$$ 
form a basis of $\mathcal V_{s+1}$. \hfill$\Box$

\begin{prop} The length generating function of the Magnus subgroup
$\mu(F_k)\subset S\mathbb K\dlangle X_1,\dots,X_k\drangle_{rat}^*$ 
is given by
$$\sum_{g\in F_k}t^{\parallel \mu(g)\parallel}=1+k(k+1)\frac{t}{1-k^2t}\ .$$

In particular, it is independent from 
$\mathbb K$ and the Magnus representation $\mu(F_k)$ 
contains exactly $(1+k)k^{2l-1}$ elements of length $l\geq 1$.
\end{prop}

{\bf Proof} By induction on $l$.
We separate elements of length $l$ in $\mu(F_k)\subset 
S\mathbb K\dlangle \mathcal X\drangle_{rat}^*$ according to 
the sign of the last exponent with respect to reduced expressions
in the free generators
$(1+X_1)^{\pm 1},\dots,(1+X_k)^{\pm 1}$. There are 
$k$ elements of length $1$ of the form $(1+X)$ and there are 
$k+k(k-1)=k^2$ elements of length $1$ of the form $1/(1+X)$ or
$(1+Y)/(1+X)$.
Let $\alpha_l$ denote the number of elements of length $l$ of the 
form $*(1+X)$. We show by induction on $l$ that $\alpha_l=k^2\alpha_{l-1}$
if $l\geq 2$ and that we have  $\beta_l=k\alpha_l$ for the number $\beta_l$
of elements of the form $*/(1+X)$ which are of length $l\geq 1$.

We have $\alpha_{l+1}=k\alpha_l+(k-1)\beta_l=(k+(k-1)k)\alpha_l=k^2\alpha_l$.
Similarly,
$$\beta_{l+1}=k\beta_l+(k-1)\alpha_{l+1}=k^2\alpha_l+(k-1)k^2\alpha_l=
k^3\alpha_l=k\alpha_{l+1}$$ 
This ends the proof.\hfill$\Box$

\begin{rem} I ignore if the metric group $\mathbb K\dlangle\mathcal X
\drangle_{rat}^*$ contains a subgroup $\mathcal G$ of finite type such
that the generating series 
$\sum_{A\in\mathcal G}t^{\parallel A\parallel}$ is irrational.

For $\mathcal G=A^{\mathbb Z}$ a non-trivial cyclic group, 
one can show rationality of the related series
$$\sum_{n=0}^\infty t^{\dim\left(\sum_{j=0}^n\overline{A^j}\right)}\ .$$
\end{rem}

\subsection{The metric monoid $S\mathbb K\dlangle \mathcal X
\drangle_{rat}^*\cap\mathbb K\langle\mathcal X\rangle$}\label{subsecmetrmon}

The set $S\mathbb K\dlangle \mathcal X\drangle_{rat}^*\cap
\mathbb K\langle\mathcal X\rangle$ is the multiplicative
monoid formed by all noncommutative polynomials with constant 
coefficient $1$.

\begin{prop} \label{propmonoidlength} (i) We have 
$$\parallel A\parallel =\dim(\overline A)-1$$
for $A\in S\mathbb K\dlangle\mathcal X\drangle_{rat}^*\cap
\mathbb K\langle \mathcal X\rangle$.

\ \ (ii) We have
$$\parallel AB\parallel=\parallel A\parallel+\parallel B\parallel$$
for $A,B\in S\mathbb K\dlangle\mathcal X\drangle_{rat}^*\cap
\mathbb K\langle \mathcal X\rangle$.
\end{prop}

\begin{lem} \label{lemmondlength} We have 
$$\overline{(AB)}=\overline A+A\overline B$$
for all $A,B\in\mathbb K\langle \mathcal X\rangle$ such that $B\not= 0$.
\end{lem}

{\bf Proof} We denote by 
$\overline{\rho(\mathcal X^{\geq n})A}$ 
the vector space generated by all series of 
the form $\rho(\mathbf X)A$ with $\mathbf X\in\mathcal X^*$ of length
$\geq n$. The vector spaces $\overline{\rho(\mathcal X^{\geq n})A}$ are
recursively closed and we have the inclusions 
$$\overline{\rho(\mathcal X^{\geq 0})A}\supset\overline{
\rho(\mathcal X^{\geq 1})A}\supset\overline{\rho(\mathcal X^{\geq 2})A}
\supset\dots\ .$$ 
A non-zero series $A\in\mathbb K\dlangle \mathcal X\drangle$ 
is a noncommutative
polynomial if and only if there exists a natural integer $D$, called
the degree of $A$, such that
$\overline{\rho(\mathcal X^{\geq D})A}=\mathbb K$ and 
$\overline{\rho(\mathcal X^{\geq D+1})A}=
\{0\}$. Assertion (i) of Lemma \ref{lemrhoXAB} implies the equalities
$$\overline{\rho(\mathcal X^{\geq n+D_{\tilde B}})(\tilde A\tilde B)}=
\overline{\rho(\mathcal X^{\geq n})\tilde A}$$
and 
$$\overline{\rho(\mathcal X^{\geq n})(\tilde A\tilde B)}=
\tilde A\overline{\rho(\mathcal X^{\geq n})\tilde B}
\pmod{\overline{\tilde A}}$$
if $\tilde B\in\mathbb K\langle \mathcal X\rangle$ is a non-zero 
polynomial of degree $D_{\tilde B}$. This proves the Lemma.\hfill$\Box$

{\bf Proof of Proposition \ref{propmonoidlength}} Assertion (i) follows from 
$\mathbb K=\overline{\rho(\mathcal X^{D_A})A}\subset 
\overline A$ and from the
definition $\parallel A\parallel=\dim(\overline A+\mathbb K)-1$.

Lemma \ref{lemmondlength} shows
$\overline{AB}=\overline A+A\overline B$.
Since $\overline A\cap A\overline B=\mathbb K A$, 
we have $\dim(\overline{AB})=\dim(\overline A)+\dim(\overline B)-1$
which shows 
$$\parallel AB\parallel=\dim(\overline A)+\dim(\overline B)-2=
\parallel A\parallel+\parallel B\parallel$$
by assertion (i).\hfill$\Box$

\subsection{The length-generating function for $\mathbb F_q\dlangle
X_1,\dots,X_k\drangle_{rat}^*$}

The aim of this Section is to give a formula 
for the generating series 
$$\sum_{A\in S\mathbb F_q\dlangle X_1,\dots,X_k\drangle_{rat}^*}
t^{\parallel A\parallel}$$
enumerating elements of $S\mathbb F_q\dlangle X_1,\dots,X_k\drangle_{rat}^*$
according to their lengths. This can be done by considering a slight variation
of the techniques and tools introduced in Section \ref{secnormal}.

As in Section \ref{secnormal}, we consider $\mathcal X^*$ as the
rooted $k-$regular infinite tree with (right-left) lexicographically
ordered vertices.

A {\it normal form} of an element 
$g\in S\mathbb K\dlangle \mathcal X\drangle_{rat}^*$
is a strictly increasing sequence $\mathbf X_1<\dots<\mathbf X_{\parallel
g\parallel}$ 
with $\mathbf X_{j+1}$ defined as the smallest element of the set
$$\{\mathbf X\in\mathcal X^*\ \vert\ \rho(\mathbf X)g\not\in
\mathbb K\oplus\bigoplus_{i=1}^j\mathbb K\rho(\mathbf X_i)g\}\ .$$ 
In particular, we have $\mathbb K+\overline g=\mathbb K\oplus
\bigoplus_{j=1}^{\parallel
g\parallel} \rho(\mathbf X_i)g$ and $\mathbf X_1=\emptyset$ if $g\not=1$.

A normal form for $g$ gives rise to a 
{\it minimal tree $G-$presentation} $(T_g,\epsilon,\mu)$ (the lettre 
$G$ stands for ``group'')
with underlying tree the finite full tree $T_g\in\mathcal{FFT}(\mathcal X)$
having interior vertices $V^\circ(T_g)$ 
given by the $\parallel g\parallel$ elements 
$\mathbf X_1,\dots,\mathbf X_{\parallel g\parallel}$ of the
normal sequence and having leaves the $1+\parallel g\parallel(k-1)$ elements
$$\partial V(T_g)=
\{\emptyset,
\mathcal X \mathbf X_1,\dots,\mathcal X\mathbf X_{\parallel g\parallel}\}
\setminus V^\circ(T_g)\ .$$
We endow interior vertices with the augmentation map 
$\epsilon:V^\circ(T_g)\longrightarrow \mathbb K$ defined by
$\epsilon(\mathbf X_j)=\epsilon(\rho(\mathbf X_j)g)\in\mathbb K$.
Since $\epsilon(g)=1$ we have always $\epsilon(\mathbf X_1)=1$
for the root vertex $\mathbf X_1=\emptyset$ of $T_g$. 
We set $\tilde V=V^\circ(T_g)\cup \{1\}$ where $\{1\}$
represents the standard basis $1$ of $\mathbb K$. The map $\mu:
\partial V(T_g)\times \tilde V\longrightarrow \mathbb K$
is defined by the equality
$$\rho(\mathbf L)g=\mu(\mathbf L,1)+\sum_{j=1}^{\parallel g\parallel}
\mu(\mathbf L,\mathbf X_j)\rho(\mathbf X_j)g\in
\mathbb K\oplus\bigoplus_{j=1}^{\parallel g\parallel} \mathbb K
\rho(\mathbf X_j)A=\mathbb K+\overline{g}\ .$$
It satisfies $\mu(\mathbf L,\mathbf Y)=0$ if $\mathbf L<\mathbf Y$ for
$\mathbf L\in\partial V(T_g)$ and $\mathbf Y\in V^\circ(T_g)$. 

Every element $g\in S\mathbb K\dlangle\mathcal X\drangle_{rat}^*$ has
a unique minimal tree $G-$presentation. A tree $T\in\mathcal{FFT}(
\mathcal X)$ underlying a tree $G-$presentation $(T,\epsilon,\mu)$
(defined in the obvious way) of an element
$g\in\mathbb K\dlangle\mathcal X\drangle_{rat}^*$ contains always
the minimal tree $T_g$ of $g$. The set of all such presentations
with underlying tree $T\in \mathcal{FFT}(\mathcal X)$
containing the minimal tree $T_g$ of $g\in \mathbb K\dlangle
\mathcal X\drangle_{rat}^*$ is an affine vectorspace of dimension
$$\sum_{\mathbf L\in \partial V(T)} \sharp\{\mathbf X\in V^\circ(T)\setminus
V^\circ(T_g)\ \vert\ \mathbf X<\mathbf L\}\ .$$

We define recursively polynomials $F_T(q)\in\mathbb N[q]$ indexed by
the set $\mathcal {FFT}(\mathcal X)$ of all full finite subtrees in  
$\mathcal X^*$ by setting
$$F_T(q)=q^{-1+\sharp(V(T))}\prod_{\mathbf L\in \partial V(T)}
q^{\sharp\{\mathbf X\in V^\circ(T)\ \vert\ \mathbf X<\mathbf L\}}-
C_T(q)$$
where 
$$C_T(q)=\sum_{T'\in\mathcal{FFT}(\mathcal X),\ T'<T}
F_{T'}(q)\prod_{\mathbf L\in\partial V(T')}q^{\sharp\{\mathbf X\in
V^\circ(T)\setminus V^\circ(T')\ \vert\ \mathbf X<\mathbf L\}}\ .$$
 
\begin{thm} We have
$$\sum_{A\in\mathbb F_q\dlangle X_1,\dots,X_k\drangle_{rat}^*}
t^{\parallel A\parallel}=\sum_{T\in\mathcal{FFT}(\mathcal X)}F_T(q)t^{
\sharp(V^\circ(T))}\ .$$

In particular, the polynomial
$$F_n(q)=\sum_{T\in\mathcal{FFT}(\mathcal X),\ \sharp(V^\circ(T))=n}
F_T(q)\in\mathbb Z[q]$$
enumerates the number of elements of length exactly $n$ in
$\mathbb F_q\dlangle X_1,\dots,X_k\drangle_{rat}^*$.
\end{thm}

\begin{exple} Working with a unique variable $X$,
one gets easily
$$\sum_{A\in\mathbb F_q[[X]]_{rat}^*}t^{\parallel A\parallel}=
\frac{1}{1-q^2t}-\frac{qt}{1-q^2t}\ .$$
In particular there exists exactly 
$q^{2n}-q^{2n-1}$ ordered pairs of polynomials $(P_1,P_2)
\in(\mathbb F_q[X])^2$ such that $P_1(0)=P_2(0)=1,\ 
\max(\deg(P_1),\deg(P_2))=n$ and $P_1,P_2$ are without common divisor.
\end{exple}

The techniques of Section \ref{sectfastcomp} can be applied if
$\mathcal X=\{X,Y\}$ and we have
$$F_n(q)=q^nw_n(q)-\sum_{j=0}^{n-1}p_{j,n}(q)F_j(q)$$
where $w_n(q)$ and $p_{j,n}(q)$ are the polynomials
defined in Section \ref{sectfastcomp}.

The first values of $F_n(q)$ are:
$$\begin{array}{l}
F_0(q)=1\\
F_1(q)=q^{2k}-q^k\\ 
F_2(q)=q^{10}+q^{9}-q^7-2q^6+q^4\\
F_3(q)=q^{18}+q^{17}+2q^{16}+q^{15}-q^{14}-2q^{13}-4q^{12}-2q^{11}
+2q^9+3q^8-q^6\\
F_4(q)={q}^{28}+{q}^{27}+2\,{q}^{26}+3\,{q}^{25}+3\,{q}^{24}+\\
\qquad +2\,{q}^{23}-{q}^{22}-4\,{q}^{21}-7\,{q}^{20}-7\,{q}^{19}
-6\,{q}^{18}-{q}^{17}+3\,{q}^{16}+\\
\qquad+5\,{q}^{15}+7\,{q}^{14}+4\,{q}^{13}+{q}^{12}-3\,{q}^{11}-4\,{q}^{
10}+{q}^{8}
\end{array}$$

For $q=2$, the first coefficients of the series 
$\sum_{n=0}^\infty F_n(2)t^n$ are 
$$\begin{array}{l}
1+12\,t+1296\,{t}^{2}+505536\,{t}^{3}+679848192\,{t}^{4}+3248147205120
\,{t}^{5}\\
+57637071142391808\,{t}^{6}+3930578658351563587584\,{t}^{7}\\+
1050888530707010579202637824\,{t}^{8}\\+
1112792971262327168651248131637248\,{t}^{9}\\+
4690276767463069086098564091958080307200\,{t}^{10}\\+
78882286441940622154458600457858710575410839552\,{t}^{11}\\+
5300169067755719965522729677599180582255569980050374656\,{t}^{12}
\end{array}$$

\begin{rem} The formulae for $E_n(q)$ and $F_n(q)$ are very similar 
and suggest to consider the common generalisation
$$P_n(q,s)=s^nw_n(q)-\sum_{j=0}^{n-1}p_{j,n}(q)P_j(q,s)\in \mathbb Z[q,s]$$
having the specialisations $E_n(q)=P_n(q,1)$ and 
$F_n(q)=P_n(q,q)$. Experimentally the specialisation
$P_n(q,1/q)$ seems to be identically $0$ for $n\geq 1$.

The specialisations $P_n(1,s)$ and $P_n(-1,s)$ have also interesting
properties.
\end{rem}

\section{A few other algebraic structures of $\mathbb K\dlangle
\mathcal X\drangle$}

This last Section surveys some related matters which are mostly
well-known, see for example \cite{BR} for a different treatment.

\subsection{Linear substitutions of variables and abelianisation} 

For a set $\mathcal X=\{X_1,\dots,X_k\}$ of $k$ variables, the 
group $\hbox{GL}_k(\mathbb K)$ of linear automorphisms of $\mathbb K^k$
acts by linear substitutions of variables on 
the algebras $\mathbb K\dlangle \mathcal X\drangle,\
\mathbb K\langle \mathcal X\rangle,\
\mathbb K\dlangle \mathcal X\drangle_{rat}$. 

\begin{prop} The natural action of $\hbox{GL}_k(\mathbb K)$ 
on $\mathbb K\dlangle\mathcal X\drangle$ by invertible linear substitutions
of the noncommutative variables $\mathcal X=\{X_1,\dots,X_k\}$
preserves the complexity.

In particular, $\hbox{GL}_k(\mathbb K)$ acts by length-preserving
automorphisms on the group $S\mathbb K\dlangle\mathcal X\drangle_{rat}^*$.
\end{prop}

We omit the easy proof.

More generally, we can substitute  the variables $X_1,\dots,X_k$ of
$\mathcal X$ by series
$M_1,\dots,M_k\in\mathfrak m$. Such a substitution
defines an endomorphism of the 
algebra $\mathbb K\dlangle\mathcal X\drangle$ which restricts 
to an endomorphism of the rational subalgebra
if and only if all series $M_1,\dots,M_k\in\mathfrak m$ are rational.

Replacing the non-commutative variables $X_1,\dots,X_k\in\mathcal X$
by commutative variables
yields  a morphism of algebras from $\mathbb K\dlangle
X_1,\dots,X_k\drangle$ onto a commutative algebra which restricts 
to a morphism from the polynomial (respectively rational) subalgebra
onto the algebra of commutative polynomials (respectively 
commutative rational fractions without singularity
at the origin). 

\begin{rem} The obvious Hankel matrix (with rows and columns indexed
by $X^\alpha Y^\beta,(\alpha,\beta)\in \mathbb N^2$) associated to the
rational 
fraction $\frac{1}{1-XY}=\sum_{n=0}^\infty X^nY^n\in\mathbb K[[X,Y]]$ 
in two commuting variables is of infinite rank. Indeed, the rational 
fractions
$\frac{1}{1-XY},\frac{X}{1-XY},\frac{X^2}{1-XY},\frac{X^3}{1-XY},
\dots$ associated to the rows $1,Y,Y^2,Y^3,\dots$ are linearly 
independent.
This behaviour is in sharp contrast with the non-commutative
case, see Example \ref{explencXY}.
\end{rem}

\subsection{The involutive antiautomorphism $\iota$}\label{sectinvolution}

Setting $\iota(X_{i_1}X_{i_2}\cdots X_{i_{l-1}}X_{i_l})=
X_{i_l}X_{i_{l-1}}\cdots X_{i_2}X_{i_1}$ for $X_{i_1}X_{i_2}\cdots 
X_{i_{l-1}}X_{i_l}\in\mathcal X^l$, the linear application
$$A=\sum_{\mathbf X\in\mathcal X^*}(A,\mathbf X)\mathbf X 
\longmapsto 
\iota(A)=\sum_{\mathbf X\in\mathcal X^*}(A,\mathbf X)\iota(\mathbf
X)$$
defines an involutive antiautomorphism $\iota$ of $\mathbb K\dlangle
X_1,\dots,X_k\drangle$ and of its polynomial and rational
subalgebras, cf. Exercise 9, Page 24 of \cite{SaSo}. 
Since the Hankel matrix $H_{\iota(A)}$ of $\iota(A)$ is essentially
the transposed matrix of the Hankel matrix 
$H_A$ of $A$, we have $\dim (\overline{\iota(A)})=\dim(\overline A)$.

Remark that the formula $\lambda(\mathbf X)A=\iota(\rho(\iota(\mathbf X))
(\iota(A)))$ defines a left-action $\lambda:\mathcal X^*\longrightarrow
\mathrm{End}(\mathbb K\dlangle\mathcal X\drangle)$ (satisfying
$\lambda(\mathbf X)(\lambda(\mathbf X')A)=\lambda(\mathbf X'\mathbf X)A$).
The dimension of the vector space spanned by the orbit 
$\lambda(\mathcal X^*)A$ equals the dimension of the space
spanned by the rows of the Hankel matrix $H_A$ for $A$ and is thus given by 
the complexity $\dim(\overline A)$ of $A$. The left and right actions
$\lambda$ and $\rho$ commute and define thus an action
$\lambda\times \rho$ of the product-monoid $\mathcal X^*\times \mathcal X^*$.
The vector space spanned by the orbit 
$(\lambda(\mathcal X^*)\times \rho(\mathcal X^*))A$ is of dimension
at most $(\dim(\overline A))^2$. More precisely, this dimension
equals the dimension of the monoid algebra 
$\mathbb K[\rho_{\overline A}(\mathcal X^*)]
\subset \mathrm{End}(\overline A)$
where $\rho_{\overline A}(\mathcal X^*)\subset \hbox{End}
(\overline A)$ denotes the shift-monoid of $A$. 

\begin{rem} One can use the left action in order to define the left-recursive
closure ${\overline A}^{\lambda}$ of an element $A\in\mathbb K\dlangle
\mathcal X\drangle$. The formula
$$A\longmapsto \parallel A\parallel_{\lambda}=\dim({\overline A}^\lambda+
\mathbb K)-1$$
turns $S\mathbb K\dlangle\mathcal X\drangle_{rat}^*$ again into a metric group
and we have 
$$\vert\ \parallel A\parallel-\parallel A\parallel_\lambda\vert\leq 1$$
for all $A\in S\mathbb K\dlangle\mathcal X\drangle_{rat}^*$. This inequality
is sharp as shown by the example 
$$\parallel (1+X)\frac{1}{(1+Y)}\parallel=1\hbox{ and }
\parallel (1+X)\frac{1}{(1+Y)}\parallel_{\lambda}=2\ .$$
\end{rem}

\subsection{Derivations}

For $X_i\in\mathcal X=\{X_1,\dots,X_k\}$ we consider the map 
$\frac{\partial}{\partial X_i}:\mathbb K\langle
\mathcal X\rangle\longrightarrow \mathbb K\langle
\mathcal X\rangle$ defined by $\frac{\partial }{\partial X_i}1=0, 
\frac{\partial }{\partial X_i}X_i=1, 
\frac{\partial }{\partial X_i}X_j=0$ for $j\not= i$ and extended
linearly to $\mathbb K\langle \mathcal X\rangle$
by the Leibnitz rule $\frac{\partial }{\partial X_i}(
\mathbf{XY})=\left(\frac{\partial}{\partial X_i}\mathbf
  X\right)\mathbf Y+
\mathbf X\left(\frac{\partial }{\partial X_i}\mathbf Y\right)$.
These maps define derivations of the polynomial algebra 
$\mathbb K\langle \mathcal X\rangle$ which extend 
in the obvious way to derivations of $\mathbb K\dlangle \mathcal X\drangle$.

A straightforward computation shows the identities
$$\rho(X)\left(\frac{\partial}{\partial Y}A\right)=
\frac{\partial}{\partial Y}\left(\rho(X)A\right)$$
and
$$\rho(X)\left(\frac{\partial}{\partial X}A\right)=
\rho(X^2)A+\frac{\partial}{\partial X}
\left(\rho(X_i)A\right)$$
for all $X,Y\in\mathcal X$ such that $X\not= Y$.
We have thus the inclusion 
$$\overline{(\partial /\partial X)A}\subset
\overline A+(\partial/\partial X)\overline{A}$$
and the resulting inequality
$$\dim(\overline{(\partial/\partial X)A})
\leq 2\dim(\overline A)$$
shows that the derivations $\partial/\partial X,\ X\in\mathcal X$, 
preserve the subalgebra $\mathbb K\dlangle \mathcal X\drangle_{rat}$
of rational elements.

\subsection{Hadamard product}\label{sechadamprod}

The Hadamard product $A\circ_H B$ of $A,B\in\mathbb K\dlangle
\mathcal X\drangle$ is defined by the coefficient-wise product
$$A\circ_HB=\sum_{\mathbf X\in\mathcal X^*}(A,\mathbf X)(B,\mathbf X)
\mathbf X\ .$$
It turns $\mathbb K\dlangle \mathcal X\drangle$ into a commutative (and
associative) algebra.

The formula $\rho(X)(A\circ_H B)=(\rho(X)A)\circ_H(\rho(X)B)$ shows
$\overline{(A\circ_H B)}\subset \overline A\circ_H\overline B$ where
$\overline A\circ_H\overline B$ denotes the vector space spanned by
Hadamard products $\tilde A\circ_H\tilde B$ with $\tilde A\in\overline A$
and $\tilde B\in\overline B$.
This implies $\dim\left(\overline{(A\circ_H B)}\right)\leq
\dim(\overline A)\dim(\overline B)$ and the Hadamard product
preserves thus the vector space $\mathbb K\dlangle\mathcal X\drangle_{rat}$
of rational elements, cf. Theorem 4.4, Page 32 of \cite{SaSo}. 
The unit group of 
$\mathbb K\dlangle\mathcal X\drangle$ for the Hadamard product
consists of all series 
involving all monomials with non-zero coefficient and identity 
the characteristic function $\sum_{\mathbf X\in \mathcal X^*}\mathbf X=
1/(1-\sum_{X\in \mathcal X}X)$ of $\mathcal X^*$.

\begin{rem} The inverse for the Hadamard product of a rational 
element involving all monomials with non-zero coefficients is in
general not rational. An example is given by 
$1+\sum_{n=1}^\infty (n+1)\left(\sum_{X\in \mathcal X}X\right)^n
\in \mathbb K\dlangle\mathcal X\drangle_{rat}$ over a field
$\mathbb K$ of characteristic $0$. 

Over the algebraically closed field $\overline{\mathbb F_p}$ of 
positive characteristic $p$ there are no such examples:
Rational series involving all monomials with non-zero coefficients
have finite order with respect to the Hadamard product.
\end{rem}

\subsection{Shuffle product}

The shuffle product is the obvious bilinear product of $\mathbb K
\dlangle\mathcal X\drangle$ defined recursively by the formulae
$1\sh \mathbf X=\mathbf X\sh 1=\mathbf X$ and 
$$(\mathbf X X)\sh(\mathbf Y Y)=\left(\mathbf X \sh (\mathbf Y Y\right)X
+\left((\mathbf X X)\sh \mathbf Y\right)Y$$
for all $\mathbf X,\mathbf Y\in\mathcal X^*$ and $X,Y\in\mathcal X$. 
The shuffle product turns the vector space $\mathbb K\dlangle 
\mathcal X\drangle$ into a commutative (and associative) algebra.

The recursive definition of the shuffle product implies 
$$\rho(X)\left(A\sh B\right)=(\rho(X)A)\sh B+A\sh (\rho(X)B)$$
which shows $\overline{(A\sh B)}\subset \overline A\sh \overline B$
where the right side denotes as usual the vector space spanned by
all elements $\tilde A\sh \tilde B$ for $\tilde A\in \overline A,
\tilde B\in \overline B$. We have thus the inequality 
$$\dim\left(\overline{A\sh B}\right)\leq \dim(\overline A)
\dim(\overline B)$$ which shows that the shuffle product restricts
to $\mathbb K\dlangle \mathcal X\drangle_{rat}$, cf. Exercice 6, Page 35 
of \cite{SaSo}. The unit group
of $\mathbb K\dlangle \mathcal X\drangle$ for the shuffle product
is the set $\mathbb K^*+\mathfrak m$ of all series with non-zero constant 
coefficient. The unit group of the rational shuffle-algebra
$\mathbb K\dlangle \mathcal X\drangle_{rat}$ is much smaller since
the shuffle inverse of a rational element in 
$\mathbb K^*+\mathfrak m\subset \mathbb K\dlangle\mathcal X\drangle_{rat}$
is in general not rational. It contains however geometric progressions
$\left(1-\sum_{j=1}^k \lambda_j X_j\right)=\sum_{n=0}^\infty \left(\sum_{
j=1}^k \lambda_j X_j\right)^n$ since we have 
$$\frac{1}{1-\sum_{j=1}^k \lambda_j X_j}\sh
\frac{1}{1-\sum_{j=1}^k \mu_j X_j}=\frac{1}{1-\sum_{j=1}^k(\lambda_j+\mu_j)
X_j}\ .$$
If the ground field $\mathbb K$ is of positive characteristic $p$, then
$$(1+a)^{\sh^p}=1$$
for $a\in\mathfrak m$
where $(1+a)^{\sh^p}$ denotes the $p-$th shuffle power 
(shuffle product $(1+a)\sh (1+a)\sh\cdots\sh (1+a)$ of
$p$ identical factors $1+a$). The shuffle inverse
of a rational element in $\mathbb K^*+\mathfrak m$ is thus again rational
in positive characteristic.

{\bf Problem} Given a field $\mathbb K$ of characteristic $0$, describe the
smallest algebra $\mathcal A\subset \mathbb K\dlangle
\mathcal X\drangle$ such that $\mathcal A$ is rationally closed
for the ordinary product and for the shuffle product.
Otherwise stated, describe the smallest algebra $\mathcal A$
which contains $\mathbb K\dlangle \mathcal X\drangle_{rat}$ such that
for every $A\in 1+\mathfrak m\cap\mathcal A$ there exist elements
$B,C\in\mathcal A$ such that $AB=1$ and $A\sh C=1$.
Remark that this algebra $\mathcal A$ is enumerable for an
enumerable field $\mathbb K$ and $\mathcal A$ is thus strictly
smaller than the non-enumerable algebra
$\mathbb K\dlangle \mathcal X\drangle$.

\subsection{Composition and homographies}

Given $A,B\in\mathbb K\dlangle \mathcal X\drangle$ where $\mathcal X=
\{X_1,\dots,X_n\}$, we set 
$$A\circ B=A(BX_1,\dots,BX_n)B$$
where $A(BX_1,\dots,BX_n)\in\mathbb K\dlangle\mathcal X\drangle$
is obtained by the substitutions $X_j\longmapsto 
BX_j,\ j=1,\dots,n$ in the non-commutative
formal power series $A$. Since $BX_j\in\mathfrak m$, 
the result of these substitutions defines a unique element
of $\mathbb K\dlangle\mathcal X\drangle$. Easy computations show that
the product defined by $(A,B)\longmapsto A\circ B$
is left-linear and associative and that it turns 
$1+\mathfrak m\subset \mathbb K\dlangle 
\mathcal X\drangle$ into a non-commutative group.

\begin{rem} The group $1+\frak m$ considered above is the diagonal 
subgroup of the group $\mathcal {FD}$ of ``formal non-commutative
diffeomorphisms tangent to the identity'' defined as follows:
$\mathcal{FD}=(1+\mathfrak m)^n\subset \mathbb K\dlangle X_1,\dots,X_n
\drangle$ as a set with product given by
$$(A_1,\dots,A_i,\dots,A_n)(B_1,\dots,B_n)=
(\dots,A_i(B_1X_1,\dots,B_nX_n)B_i,\dots)\ .$$
The group law on $\mathcal{FD}$ is composition where 
$(A_1,\dots,A_n)$ corresponds to the formal diffeomorphism 
$(A_1X_1,\dots,A_nX_n)$.

One could of course also consider compositions of elements
of the form $(X_1A_1,\dots,X_nA_n)$. The resulting group is 
isomorphic to $\mathcal{FD}$.
\end{rem}

The formula 
$$\rho(X)(A\circ B)=A(BX_1,\dots,BX_n)\left(\rho(X)B\right)+
\epsilon(B)\left(\rho(X)A\right)(BX_1,\dots,BX_n)B$$
and left linearity of the compositional product $(A,B)\longmapsto
A\circ B$
show that 
$$\overline{A\circ B}\subset \overline A(BX_1,\dots,BX_n)\overline B$$
which implies $\dim(\overline{A\circ B})\leq\dim(\overline A)
\dim(\overline B)$. The compositional product $\circ$ turns thus
the set $\mathbb K\dlangle \mathcal X\drangle_{rat}$ of rational elements 
and its subset $1+\mathfrak m\cap \mathbb K\dlangle\mathcal X\drangle_{rat}$
into monoids. The compositional inverse $B$, defined by $B\circ A=A\circ B=1$,
of a rational element $A\in 1+\mathfrak m\cap \mathbb K\dlangle 
\mathcal X\drangle_{rat}$ is in general not rational.
However, a straightforward computation yields
$$\begin{array}{l}
\displaystyle \frac{1}{1-\sum_{j=1}^k\lambda_jX_j}\circ\frac{1}{1-\sum_{j=1}^k\mu_jX_j}\\
\displaystyle \quad =
(1-\sum_{j=1}^k\mu_jX_j)\frac{1}{1-\sum_{j=1}^k(\lambda_j+\mu_j)X_j}\frac{1}
{1-\sum_{j=1}^k\mu_jX_j}\end{array}$$
and shows that the compositional inverse of a geometric 
progression given by the rational series 
$$\frac{1}{1-\sum_{j=1}^k\lambda_jX_j}=1+\sum_{n=1}^\infty\left(\sum_{j=1}^k
\lambda_jX_j\right)^n$$ is the rational series $1/(1+\sum_{j=1}^k
\lambda_jX_j)$. We call the subgroup $\mathcal H$ of the compositional 
group $1+\mathfrak m$ generated by all rational elements of the 
form $1/(1-\sum_{j=1}^k
\lambda_jX_j),\ (\lambda_1,\dots,\lambda_k)\in\mathbb K^k$, the 
{\it group of homographies}. It would be interesting to know if there
exist rational elements $A,B\in
1+\mathfrak m\cap \big(\mathbb K\dlangle\mathcal X\drangle_{rat}\setminus 
\mathcal H\big)$ such that $A\circ B=1$.

{\bf Problem} As for the shuffle product, one might ask to describe the 
smallest subalgebra $\mathcal A\subset \mathbb K\dlangle \mathcal X\drangle$
which contains $\mathbb K\dlangle X\drangle_{rat}$ and
intersects $1+\mathfrak m$ in a subgroup for the compositional product.
One might in fact ask for characterising the smallest rationally closed
algebras which are ``closed'' with respect to the corresponding 
group structure given by one or more of the monoid structures 
associated to the Hadamard product, the shuffle product and the compositional
product. The largest such algebra, defined as being closed with respect 
to inversion of invertible elements for all 
four group-laws (ordinary non-commutative product, Hadamard product,
shuffle product and compositional product) is enumerable over 
an enumerable field $\mathbb K$ and thus distinct from 
$\mathbb K\dlangle \mathcal X\drangle$.

\subsection{Automatic sequences}

This section gives a very brief outline without details or proofs
of the link between certain 
rational elements of $\mathbb K\dlangle\mathcal X\drangle$
and so-called automatic sequences, see \cite{AS} for the definition.

Given a natural integer $k\geq 2$, we can consider the injection
$\mathbb K^{\mathbb N}\longrightarrow\mathbb K\dlangle X_0,\dots,X_{k-1}
\drangle$ given by the map
$$(s(0),s(1),\dots)\longmapsto
\sum_{\mathbf X=X_{i_0}\dots X_{i_l}\in\{X_0,\dots,X_{k-1}\}^*}
s\left(\sum_{j=0}^l i_jk^j\right)\mathbf X$$
or the bijection $\mathbb K^{\mathbb N}\longrightarrow\mathbb K\dlangle
X_1,\dots,X_k\drangle$ defined by
$$(s(0),s(1),\dots)\longmapsto
\sum_{\mathbf X=X_{i_0}\dots X_{i_l}\in\{X_1,\dots,X_k\}^*}
s\left(\sum_{j=0}^l i_jk^j\right)\mathbf X\ .$$
Let $\mathcal I\subset \mathbb K\dlangle\mathcal X\drangle$ denote
the image of one of these maps.
The image in $\mathcal I$ of the set of $k-$automatic sequences in
$\mathbb K^{\mathbb N}$ is then exactly the subset $\mathcal I_f
\subset \mathcal I\cap\mathbb K\dlangle\mathcal X\drangle_{rat}$ 
corresponding to rational elements
with coefficients in a finite subset of the field
$\mathbb K$. By Proposition \ref{propfinmonoid},
a rational element $A\in\mathbb K\dlangle\mathcal X\drangle_{rat}$
has its coefficients in a finite subset of $\mathbb K$ if and only if 
it has a finite shift monoid $\rho_{\overline A}(\mathcal X^*)\subset
\hbox{End}(\overline A)$. Using correct conventions,
a finite-state automaton
for the $k-$automatic sequence associated to such an element $A\in I_f\cap
\mathbb K\dlangle\mathcal X\drangle_{rat}$ is given by the Cayley 
graph (with respect to the generators $\rho_{\overline A}
(\mathcal X)$) of the
finite monoid $\rho_{\overline A}(\mathcal X^*)\subset 
\hbox{End}(\overline A)$.
The initial state of the finite state automaton is $\rho_{\overline A}
(\emptyset)$ and the output function 
$\rho_{\overline A}(\mathbf X)\longmapsto \epsilon(\rho_{
\overline A}(\mathbf X)A)$, 
see Chapter 4 of \cite{AS} for definitions.

\subsection{Regular languages}

A {\it language} is a subset of $\mathcal X^*$ over a finite alphabet
$\mathcal X$.

A {\it finite-state automaton} is a finite oriented graph $\Gamma$
such that:

$\Gamma$ contains a marked initial vertex $v_*$.

Each vertex of $\Gamma$ is the initial vertex of exactly 
$\sharp(\mathcal X)$ oriented edges, labelled by $\mathcal X$.

The vertices of $\Gamma$ are partitioned into two finite disjoint
subsets $\mathcal A$ and $\mathcal R$.

A finite-state automaton $\Gamma$ defines a unique language 
$\mathcal L(\Gamma)$,
called the language accepted by $\Gamma$, as follows:
Every word $X_{i_1}\dots X_{i_l}$ of $\mathcal X^*$
defines a unique oriented path
starting at $v_*$ and consisting of the $l$ consecutive oriented 
edges labelled $X_{i_l}, X_{i_{l-1}},\dots,X_{i_2},X_{i_1}$.
The word $X_{i_1}\dots X_{i_l}\in\mathcal X^*$ belongs to 
$\mathcal L(\Gamma)$ is
and only if the associated path ends in a vertex of the subset 
$\mathcal A$ of {\it accepting states}. 

A language $\mathcal L\subset\mathcal X^*$ is {\it regular} (some 
authors say also {\it rational} or {\it recognisable}, cf \cite{BR})
if it is accepted by a finite-state automaton.

A recursive presentation $A_j=\gamma_j+\sum_j A_i\alpha_{i,j},\ j
\in\mathcal I$ (with $\mathcal I$ finite), of a series $A=A_1
\in \mathbb R\dlangle\mathcal X\drangle$
is {\it positive} if $\gamma_j\geq 0,\ \rho(X)\alpha_{i,j}\geq 0$ for all 
$i,j\in\mathcal I$ and for all $X\in\mathcal X$.

Such a recursive presentation is {\it integral} if 
$\gamma_j\in\mathbb Z,\ \rho(X)\alpha_{i,j}\in \mathbb Z$ for all 
$i,j\in\mathcal I$ and for all $X\in\mathcal X$.

The following result is also contained in Chapter III of \cite{BR}
or in Section II.5 of \cite{SaSo}:

\begin{prop} \label{propcharreglang}
The following statements are equivalent:

\ \ (i) $\mathcal L\subset \mathcal X^*$ is a regular language.

\ \ (ii) The characteristic function
$$\sum_{\mathbf X\in\mathcal L}\mathbf X$$
of $\mathcal L\subset\mathcal X^*$
is a rational series of $\mathbb K\dlangle\mathcal X\drangle$
for any field $\mathbb K$.

\ \ (iii) The characteristic function
$$\sum_{\mathbf X\in\mathcal L}\mathbf X$$
of $\mathcal L\subset\mathcal X^*$
is a rational series of $\mathbb K\dlangle\mathcal X\drangle$
for some field $\mathbb K$.

\ \ (iv) $\mathcal L$ is the support of a rational series in 
$\mathbb Q\dlangle
\mathcal X\drangle$ which has an integral positive recursive presentation.

\ \ (v) $\mathcal L$ is the support of a rational series in $\mathbb R
\dlangle \mathcal X\drangle$
having a positive recursive presentation.
\end{prop}

Since rational series are closed under Hadamard products
(see Section \ref{sechadamprod}),
we have:

\begin{cor} The set of all regular languages is also closed under
intersections and differences.
\end{cor}

\begin{rem} The positivity conditions in assertions (iv) and (v) are
necessary as shown by examples in \cite{BR}.
\end{rem}

{\bf Proof of Proposition \ref{propcharreglang}} A finite-state 
automaton $\Gamma$ for a regular language $\mathcal L$ 
defines a recursive presentation for $\mathcal L$ as follows:
consider the series $A_v\in\mathbb K\dlangle\mathcal X\drangle$ 
indexed by vertices $v\in V(\Gamma)$ of $\Gamma$
which are defined by the equations
$$A_v=\epsilon(v)+\sum_{\mathbf X\in\mathcal X}\rho(X)A_v,\ 
v\in V(\Gamma)$$
where $\epsilon(v)=1$ if $v\in\mathcal A$ and $\epsilon(v)=0$ otherwise
and where $\rho(X)A_v=A_w$ if an oriented edge labelled $X$ starts
at $v$ and ends at $w$. We have then clearly 
$A_{v_*}=\sum_{\mathbf X\in\mathcal L}\mathbf X$.
This shows that (i) implies (ii).

Assertion (ii) implies (iii) trivially.

Consider a rational function of the form
$A=\sum_{\mathbf X\in\mathcal L}\mathbf X$ for  $\mathcal L\subset
\mathcal X^*$. Proposition \ref{propfinmonoid} shows that 
$\rho_{\overline A}(\mathcal X^*)$ is finite.
The finite set $\rho_{\overline A}(\mathcal X^*)A$ is thus stable 
under shift-maps 
and since $\epsilon(\rho(\mathcal X^*)A)\subset \{0,1\}$,
it can be used for writing down a presentation using only coefficients
in $\{0,1\}\subset \mathbb R_{\geq 0}$. This shows that
(iii) implies (iv). 

Assertion (iv) implies obviously (v).

Given a presentation of $A$ involving only non-negative real numbers,
we can use the Boolean algebra $\mathbb B=\{0,p\}$ defined by 
$0+0=0,0+p=p+0=p+p=p$ and $0\cdot 0=0\cdot p=p\cdot 0=0,\ 
p\cdot p=p$ in order to define an element
$\tilde A\in\mathbb B\dlangle \mathcal X\drangle$
which has the same support as $A$ by replacing each strictly
positive real number arising in the recursive presentation of $A$ by $p$.

The resulting shift monoid $\tilde\rho(\mathcal X^*)$
over the algebra $\mathbb B$ is finite.
The finite state automaton given by its Cayley graph with 
accepting states $\mathcal A$ defined by $\tilde\rho(\mathbf X)\in
\mathcal A$ if $\epsilon(\rho(X)A)>0$ is a finite state automaton
with accepted language the support of $A$.
This shows that (v) implies (i) and ends the proof.\hfill$\Box$

We end this brief section by mentionning a last well-known result:

\begin{prop} The set of all regular languages is the smallest subset of
$\mathcal P(\mathcal X^*)$ which contains all finite subsets and which is
closed under unions, concatenations and the Kleene closure 
$\mathcal L\longmapsto \mathcal L^*=\cup_{n=0}^\infty \mathcal L^n$.
\end{prop}

{\bf Proof} If $\mathcal L,\mathcal L'\subset \mathcal X^*$
are two regular languages given as supports of rational
series $A,A'\in\mathbb R\dlangle\mathcal X\drangle$ having positive
presentations with respect to finite sets $A_1=A,A_2,\dots$
and $A'_1=A',A'_2,\dots$ spanning $\overline A$ and $\overline{A'}$,
then $\mathcal L\cup \mathcal L'$, respectively 
$\mathcal L\mathcal L'$, is the support of $A+A'$, respectively $AA'$,
having a positive presentation with respect to 
$A+A',A_1,A_2,\dots,A'_1,A'_2,\dots$, respectively $A_iA'_j$.
If $\mathcal L$ is regular then $\tilde{\mathcal L}=
\mathcal L\setminus \{\emptyset\}$
is also regular and $\mathcal L^*=\tilde{\mathcal L}^*$. 
We suppose thus $\emptyset\not\in \mathcal L$
and consider the characteristic function $A=\sum_{\mathbf X\in\mathcal L}
\mathbf X$ having by Proposition \ref{propfinmonoid}
a finite orbit $\rho(\mathcal X^*)A$.
The formula given by assertion (ii) of Lemma \ref{lemrhoXAB}
implies then that $B=1/(1-A)$ has a positive presentation
with respect to the finite set $B,B\rho(\mathcal X^*)A$
and the support of $B$ is obviously the Kleene closure of $\mathcal L$.

This shows that every language obtained by unions, concatenations
and Kleene closures from finite subsets in $\mathcal X^*$ is regular.

The opposite direction is given by inspecting the proof of Proposition
\ref{proprecognimpliesrational}, applied to a positive 
presentation.\hfill$\Box$

{\bf Acknowledgements} I thank P. de la Harpe for his interest and comments.
I thank also C. Reutenauer very strongly for pointing out many 
inaccuracies and omissions in a first version.


\noindent Roland BACHER

\noindent INSTITUT FOURIER

\noindent Laboratoire de Math\'ematiques

\noindent UMR 5582 (UJF-CNRS)

\noindent BP 74

\noindent 38402 St Martin d'H\`eres Cedex (France)
\medskip

\noindent e-mail: Roland.Bacher@ujf-grenoble.fr

\end{document}